\documentclass[reqno,11pt]{amsart}

\usepackage{epsf,amssymb,graphics,graphicx,wrapfig}

\newcommand{\floor}[1]{\left\lfloor #1 \right\rfloor}
\def\eps{\varepsilon}

\def\N{\bbN}

\def\be{\begin{equation}}
\def\ee{\end{equation}}
\def\ba{\begin{align}}
\def\bm{\begin{multline}}
\def\bfig{\begin{figure}[htb]}
\def\efig{\end{figure}}

\numberwithin{equation}{section}
\newtheorem{theorem}{Theorem}[section]
\newtheorem{proposition}[theorem]{Proposition}
\newtheorem{lemma}[theorem]{Lemma}
\newtheorem{corollary}[theorem]{Corollary}

\DeclareMathSymbol{\leqslant}{\mathalpha}{AMSa}{"36}
\DeclareMathSymbol{\geqslant}{\mathalpha}{AMSa}{"3E}
\DeclareMathSymbol{\doteqdot}{\mathalpha}{AMSa}{"2B}
\DeclareMathSymbol{\circlearrowright}{\mathalpha}{AMSa}{"08}
\DeclareMathSymbol{\subsetneq}{\mathalpha}{AMSb}{"28}
\DeclareMathSymbol{\supsetneq}{\mathalpha}{AMSb}{"29}
\renewcommand{\leq}{\;\leqslant\;}
\renewcommand{\geq}{\;\geqslant\;}
\newcommand{\isdefby}{\;\doteqdot\;}

\newcommand{\dd}{{\rm d}}
\newcommand{\e}[1]{\,{\rm e}^{#1}\,}
\newcommand{\ii}{{\rm i}}
\newcommand{\sumtwo}[2]{\sum_{\substack{#1 \\ #2}}}
\newcommand{\sumthree}[3]{\sum_{\substack{#1 \\ #2 \\ #3}}}

\def\dist{{\operatorname{dist\,}}}
\def\diam{{\operatorname{diam}}}
\def\Re{{\operatorname{Re\,}}}

\newcommand{\upchi}{\raise 2pt \hbox{$\chi$}}
\newcommand{\smallupchi}{\raise 1pt \hbox{\tiny $\chi$}}

\def\writefig#1 #2 #3 {\rlap{\kern #1 truecm \raise #2 truecm
\hbox{#3}}}

\newcommand{\caN}{{\mathcal N}}\newcommand{\caS}{{\mathcal S}}

\newcommand{\bbN}{{\mathbb N}}\newcommand{\bbR}{{\mathbb R}}\newcommand{\bbZ}{{\mathbb Z}}

\newcommand{\bsk}{{\boldsymbol k}}\newcommand{\bsn}{{\boldsymbol n}}\newcommand{\bsx}{{\boldsymbol x}}\newcommand{\bsy}{{\boldsymbol y}}

\newcommand{\bsN}{{\boldsymbol N}}

\newcommand{\bsvarrho}{{\boldsymbol\varrho}}

\begin{document}


\title[Spatial random permutations with small cycle weights]{Spatial random permutations \\ with small cycle weights}

\author{Volker Betz and Daniel Ueltschi}
\address{Volker Betz and Daniel Ueltschi \hfill\newline
\indent Department of Mathematics \hfill\newline
\indent University of Warwick \hfill\newline
\indent Coventry, CV4 7AL, England \hfill\newline
{\small\rm\indent http://www.maths.warwick.ac.uk/$\sim$betz/} \hfill\newline
{\small\rm\indent http://www.ueltschi.org}
}
\email{v.m.betz@warwick.ac.uk, daniel@ueltschi.org}

\maketitle

\begin{quote}
{\small
{\bf Abstract.}
We consider the distribution of cycles in two models of random permutations, that are related to one another.
In the first model, cycles receive a weight that depends on their length.
The second model deals with permutations of points in the space and there is an additional weight that involves the length of permutation jumps.
We prove the occurrence of infinite macroscopic cycles above a certain critical density.
}  

\vspace{1mm}

\vspace{1mm}
\noindent
{\footnotesize {\it 2000 Math.\ Subj.\ Class.:} 60K35, 82B20, 82B26, 82B41.}
\end{quote}

\section{Introduction}

Random permutations and their cycle structure have been studied for many decades, with a strongly increased activity in recent years. Apart from the rich mathematical structure, this interest is justified by a wide range of applications, 
from Gromov-Witten theory \cite{Ok2005} to polynuclear growth \cite{FSP} and to mathematical biology \cite{Ewe}. Detailed properties have been established in the case of uniform permutations for the moments of the distribution of the $n$-th shortest (or longest) cycle \cite{ShLo66}, for the longest increasing subsequence \cite{BDJ}, or for the convergence to equilibrium \cite{Sch}. Physics and biology have suggested certain models with nonuniform permutations for which interesting results have been obtained \cite{Pit,FH}.

In the present work, we introduce a model for nonuniform permutations which is motivated by its connection to the theory of Bose-Einstein condensation \cite{Fey,Suto1,GRU,BU}. Mathematically, the distinguishing feature of our model, when compared to the works mentioned above, is that the measure on permutations
possesses a spatial structure. More precisely, we consider pairs $(\bsx,\pi)$ with $\bsx \in \Lambda^N$ ($\Lambda$ is a cubic box in $\bbR^d$) and $\pi \in \caS_N$ (the group of permutations of $N$ elements). The weight of $(\bsx,\pi)$ is given by the ``Gibbs factor'' $\e{-H(\bsx,\pi)}$ with Hamiltonian of the form
\be \label{energy}
H(\bsx,\pi) = \sum_{i=1}^N \xi(x_i-x_{\pi(i)}) + \sum_{\ell\geq1} \alpha_\ell r_\ell(\pi).
\ee
We always assume that $\xi$ is a function $\bbR^d \to \bbR \cup \{\infty\}$, with $\int\e{-\xi(x)} \dd^d x = 1$. 
The cycle parameters $\alpha_1,\alpha_2,\dots$ are some fixed numbers, typically but not necessarily positive, and $r_\ell(\pi)$ is the number of cycles of length $\ell$ in the permutation $\pi$. The length of the cycle that contains the index $i$ is the smallest integer $\ell\geq1$ such that $\pi^\ell(i)=i$ --- this definition of the length involves the permutation but not the underlying spatial structure. Intuitively, the Gibbs factor restricts the permutations so each jump is local, i.e.\ the distances $|x_i - x_{\pi(i)}|$ remain finite even for large systems. The main question deals with the lengths of the cycles, in the limit of infinite volumes.

When $\xi(x) = \gamma |x|^2 + c$ and $\alpha_\ell = 0$ for all $\ell$, we obtain the model of spatial random permutations that corresponds to the ideal Bose gas; in this case $\rho_{\rm c}$ is the well-known critical density for Bose-Einstein condensation for non-interacting particles. The occurrence of macroscopic cycles in the ideal Bose gas has been understood in \cite{Suto1,Suto2}. The present setting with general functions $\xi$ was considered in
\cite{BU}. The latter article also introduces the cycle weight $\alpha_2$ as an approximation for the interactions between quantum particles, and the occurrence of macroscopic cycles was proved for large densities. The present article extends the results of \cite{BU} to more general cycle weights and to all densities larger than the critical density.

The main result of this article deals with the occurrence of infinite permutation cycles. Let $N_k(\pi)$ be the random variable that counts the number of points in cycles of length $k$ (we have $\sum_k N_k(\pi) = N$ for all $\pi$), and let $E_{\Lambda,N}(N_k)$ be its expectation. We consider the thermodynamic limit $|\Lambda|,N \to \infty$ with fixed density $\rho = N/|\Lambda|$ ($|\Lambda|$ denotes the volume of $\Lambda$). Fatou's lemma implies that
\be
\sum_{k\geq1} \lim_{|\Lambda|\to\infty} E_{\Lambda,\rho |\Lambda|} \Bigl( \frac{N_k}{|\Lambda|} \Bigr) \leq \rho.
\ee
We prove in this article that the left side is {\it strictly less} than $\rho$ if and only if the density is larger than a {\it critical density} $\rho_{\rm c} \in (0,\infty]$. The precise formulation of this result can be found in Theorem \ref{thminfinitecycles}. We need to restrict to certain functions $\xi$ (namely, $\e{-\xi}$ has positive Fourier transform) and small cycle weights, in the sense that $\alpha_\ell \to 0$ as $\ell \to \infty$, faster than $1/\log\ell$.
Our results include an explicit formula for $\rho_{\rm c}$, cf.\ Eq.\ (\ref{critdens}), and some characterization of the nature of infinite cycles --- they are macroscopic and the distribution of a given macroscopic cycle is uniform in $[0,\rho-\rho_{\rm c}]$. Points that are not in macroscopic cycles are shown to be necessarily in finite cycles, and their density is given by $\max(\rho,\rho_{\rm c})$.
To our knowledge, the presence of a density of points in finite cycles is specific to spatial models, and does not 
occur in the other known models of random permutations.

The model considered here can be thought of as the ``annealed" version of another ``quenched" model where the space positions would be fixed. The latter model seems more attractive from a probability perspective, and it is actually discussed in \cite{BU}. There are two main advantages of the present annealed model: the availability of rigorous results about infinite cycles, and its closeness to the Feynman-Kac representation of the quantum Bose gas.

The structure of this article is as follows. In Section \ref{sec cycle weight} we introduce an auxiliary model of 
non-spatial permutations with cycle weights, which turns out to be closely related to the spatial one. It corresponds to taking $\bsx=(0,0,\dots)$. We only consider the case of small weights $\alpha_\ell$ so that the typical nonspatial permutations are like those with uniform distribution. We will study more general weights and other behaviors in a subsequent article \cite{BUV}.

In Section \ref{sec spatial model} we introduce the model of spatial permutations with cycle weights. We discuss the existence of the infinite volume limits for thermodynamic potentials and the equivalence of ensembles. These notions belong to statistical mechanics rather than probability theory, but we need these results when we consider the more relevant question, as far as probability theory is concerned, of the occurrence of infinite cycles. Our main result is Theorem \ref{thminfinitecycles} in Section \ref{prob model}. We find in particular that the cycle weights modify the critical density; so they do have an effect on spatial permutations, unlike what was observed in Section \ref{sec cycle weight}.

Section \ref{sec thermo limit} is devoted to the proofs of the results about the thermodynamic potentials. We adapt the classical methods of Fisher and Ruelle \cite{Rue} to our context. In Section \ref{sec Fourier} we relate our spatial model to an equivalent model in the Fourier representation, which allows to prove our main theorem. The nonspatial model with cycle weights also plays a r\^ole here. We use techniques introduced for the ideal Bose gas in \cite{BP} and \cite{Suto2}.

\section{The simple model of random permutations with cycle weight}
\label{sec cycle weight}

We start with the study of random permutations of $n$ elements with no spatial structure, but with cycle weights. The results of this section will be useful for the spatial model and this is our main motivation. But the present model has its own interest. We consider only small cycle weights here, but results for other regimes will be presented in a subsequent article \cite{BUV}.

\subsection{Setting and properties}

The sample space is $\caS_n$ and the probability of a permutation $\pi \in \caS_n$ is given by
\be
p_n(\pi) = \frac1{h_n n!} \exp\Bigl\{ -\sum_{\ell\geq1} \alpha_\ell r_\ell(\pi) \Bigr\}
\ee
with normalization
\be
\label{def h_n}
h_n = \frac1{n!} \sum_{\pi \in \caS_n} \e{-\sum_\ell \alpha_\ell r_\ell(\pi)}.
\ee
Here,  $\alpha_1, \alpha_2, \dots$ are fixed numbers and $r_\ell(\pi)$ denotes the number of cycles of length $\ell$ in the permutation $\pi$. Notice the symmetry: since $\sum_j j r_j(\pi) = n$ for all $\pi$, the probability $p_n$ is invariant under the transformation
\be
\alpha_j \mapsto \alpha_j + cj,
\ee
for any constant $c$; the normalization satisfies $h_n \mapsto \e{-cn} h_n$.

Let $N_{a,b}(\pi) = \sum_{\ell=a}^b \ell r_\ell(\pi)$ denote the number of indices that belong to cycles of length between $a$ and $b$. Our main result deals with the asymptotic distribution of cycle lengths.

\begin{theorem}
\label{thm cycle weight}
If $\sum_{\ell\geq1} \frac1\ell |1-\e{-\alpha_\ell}| < \infty$, we have for any $0 \leq s \leq 1$
\[
\lim_{n\to\infty} \tfrac1n E_n(N_{1,sn}) = s.
\]
\end{theorem}

In essence, the hypothesis of the theorem requires that $\alpha_\ell \to 0$ a bit faster than $1/\log\ell$. Theorem \ref{thm cycle weight} implies that almost all indices belong to cycles whose length is a positive fraction of $n$. It can be shown that the number of cycles is of order $\log n$. The claim is easy to get in the case of uniform random permutations ($\alpha_\ell \equiv 0$), but the extension to even small weights requires some efforts. The key to the proof of Theorem \ref{thm cycle weight}, which is given in the next subsection, is the following relation. 

\begin{lemma}
\label{lem a useful relation}
\[
\label{a useful relation indeed}
E_n(N_{a,b}) = \sum_{j=a}^b \e{-\alpha_j} \frac{h_{n-j}}{h_n}.
\]
\end{lemma}

For the proof of this lemma, we just remark that $E_n(N_{a,b}) = n \, p_n(\ell_1 \in [a,b])$, with $\ell_1 = \ell_1(\pi)$ the length of the cycle that contains 1. Summing over all possible values $j$ of $\ell_1$, and observing that there are $\frac{(n-1)!}{(n-j)!}$ possible cycles, we get the relation above.

\subsection{Properties of the normalization $h_n$}

In view of Lemma \ref{lem a useful relation} it is clear that we need to gather some information on $h_n$. 
We start with a few exact relations (Proposition \ref{prop h_n formula}) and we then obtain estimates (Proposition \ref{prop estimates h_n}).

\begin{proposition}
\label{prop h_n formula}
The $h_n$'s satisfy the following properties
\begin{itemize}
\item[(a)] A recursion formula:
\[
h_n = \frac1n \sum_{\ell=1}^n \e{-\alpha_\ell} h_{n-\ell}, \qquad h_0 = 1.
\]
\item[(b)] An explicit formula:
\[
h_n = \sum_{k=1}^n \frac{1}{k!} \sumtwo{\ell_1, \ldots, \ell_k \geq 1}{\ell_1 + \ldots + \ell_k = n} \prod_{i=1}^k \frac{\e{-\alpha_{\ell_i}}}{\ell_i}, \qquad n\geq1.
\]
\item[(c)] The increments satisfy
\[
h_n - h_{n-1} = \sum_{k=1}^n \frac{1}{k!} \sumtwo{\ell_1, \ldots, \ell_k \geq 1}{\ell_1 + \ldots + \ell_k = n} \prod_{i=1}^k \frac{\e{-\alpha_{\ell_i}} - 1}{\ell_i}.
\]
\item[(d)] Another formula for $h_n$:
\[
h_n = \sum_{k=0}^n \frac1{k!} \sumtwo{\ell_1,\dots,\ell_k\geq1}{\ell_1+\ldots+\ell_k \leq n} \prod_{i=1}^k \frac{\e{-\alpha_{\ell_i}} - 1}{\ell_i}.
\]
\item[(e)] If $\frac{\alpha_\ell}\ell \to 0$ and $\gamma>0$, we have
\[
\sum_{n\geq0} \e{- \gamma n} h_n = \exp \sum_{j \geq 1} \frac{\e{-\gamma j - \alpha_j}}{j}.
\]
\end{itemize}
\end{proposition}

The formula (b) shows that $h_n$ is decreasing with respect to $(\alpha_\ell)$. The formula (d) shows that $h_n$ is increasing with respect to $n$ if $\alpha_\ell \leq 0$. Now we cannot resist but ask the reader to consider the following expression:
\[
\sum_{k=1}^n \frac1{k!} \sumtwo{\ell_1,\dots,\ell_k \geq 1}{\ell_1+\dots+\ell_k=n} \frac1{\ell_1 \dots \ell_k}.
\]
How does it behave for large $n$? The answer is surprisingly simple and is given by Proposition \ref{prop h_n formula} (b).

\begin{proof}
The recursion formula is obtained from Lemma \ref{lem a useful relation} by noting that $E_n(N_{1,n}) = n$.
For the claim (b), it is useful to define $b_j = \e{-\alpha_{j+1}}$, $j \geq 0$. The recursion formula can be written as
\be
(n+1) h_{n+1} = \sum_{j=0}^{n} b_j h_{n-j}.
\ee
Thus the series $((n+1) h_{n+1})$ is equal to the convolution of the series $(b_n)$ and $(h_n)$.
We introduce the generating functions
\be
G_h(s) = \sum_{n\geq0} h_n s^n, \qquad G_b(s) = \sum_{n\geq0} b_n s^n.
\ee
The generating function for the series $((n+1) h_{n+1})$ is $G_h'(s)$. By the properties of convolutions, we have
\be
\label{ode for G_h}
G_h'(s) = G_b(s) G_h(s), \qquad G_h(0) = 1.
\ee
The solution is
\be \label{genfunct}
G_h(s) = \exp \int_0^s G_b(t) \dd t = \exp \sum_{n\geq1} \tfrac{b_{n-1}}n s^n.
\ee
Expanding the exponential and rearranging the terms allows to find an expression for each coefficient $h_n$. This gives Proposition \ref{prop h_n formula} (b).

The generating function for the increments, $G_{\delta h}$, satisfies
\be
G_{\delta h}(s) \isdefby \sum_{n\geq0} (h_n - h_{n-1}) s^n = (1-s) G_h(s)
\ee
(with $h_{-1} = 0$). Using \eqref{ode for G_h}, we get a differential equation for $G_{\delta h}$, namely
\be
G_{\delta h}'(s) = G_{\delta h}(s) \bigl[ G_b(s) - \tfrac1{1-s} \bigr]
\ee
The expression in the bracket is equal to the generating function $G_{b-1}$ of the series $(b_n-1)$. Solving the differential equation, we get
\be
G_{\delta h}(s) = \exp \int_0^s G_{b-1}(t) \dd t = \exp \sum_{n\geq1} \frac{b_{n-1} - 1}n s^n.
\ee
Again expanding the exponential and matching the coefficients, we get Proposition \ref{prop h_n formula} (c).
The formula (d) follows from (c) since $h_n = \sum_{i=0}^n (h_i-h_{i-1})$. Finally, (e) follows directly from 
\eqref{genfunct}.
\end{proof}

We now collect a few estimates for $h_n$.

\begin{proposition}\hfill
\label{prop estimates h_n}
\begin{itemize}
\item[(a)] If $\sum_{\ell\geq1} \frac1\ell |1-\e{-\alpha_\ell}| < \infty$, we have $\sum_n |h_n-h_{n-1}| < \infty$, and
\[
h_\infty \doteqdot \lim_{n\to\infty} h_n = \exp \sum_{\ell\geq1} \frac{\e{-\alpha_\ell} - 1}\ell.
\]
\item[(b)] $h_n \geq \frac1n h_{n-1} \e{-\alpha_1}$; iterating, $h_n \geq \frac1{n!} \e{-n\alpha_1}$.
\item[(c)] $h_n \geq \frac1n \e{-\alpha_n}$.
\item[(d)] If $(\alpha_\ell)$ is subadditive (i.e.\ $\alpha_{\sum_i \ell_i} \leq \sum_i \alpha_{\ell_i}$), $h_n \leq \e{-\alpha_n}$.
\item[(e)] If $(\alpha_\ell)$ is superadditive, $h_n \geq \e{-\alpha_n}$.
\end{itemize}
\end{proposition}

\begin{proof}
The claim (a) is an immediate consequence of Proposition \ref{prop h_n formula} (d) and the dominated convergence theorem: In the limit $n \to \infty$, the constraint $\ell_1 + \ldots \ell_k \leq n$ vanishes 
and the corresponding expression factorizes. (b) and (c) follow from the recursion formula of Proposition \ref{prop h_n formula} (a), keeping only the term $\ell=1$ for (b), and the term $\ell=n$ for (c). For (d) and (e), we write the formula of Proposition \ref{prop h_n formula} (b) as
\be
h_n = \sum_{k=1}^n \frac1{k!} \sumtwo{\ell_1,\dots,\ell_k \geq 1}{\ell_1+\dots+\ell_k=n} \frac{\e{-\sum_{i=1}^k \alpha_{\ell_i}}}{\ell_1 \dots \ell_k}.
\ee
We replace $\sum\alpha_{\ell_i}$ by $\alpha_n$, getting an upper bound for $h_n$ if $(\alpha_\ell)$ is subadditive, and a lower bound if $(\alpha_\ell)$ is superadditive.
\end{proof}

\begin{proof}[Proof of Theorem \ref{thm cycle weight}]
We split the expression of Lemma \ref{lem a useful relation} into
\be
\label{a useful expanded relation indeed}
\frac1n E(N_{1,sn}) = s + \sum_{j=1}^{sn} \frac1n \frac{h_{n-j}}{h_n} (\e{-\alpha_j}-1) + \frac1n \sum_{j=1}^{sn} \Bigl( \frac{h_{n-j}}{h_n} - 1 \Bigr).
\ee
We suppose that $sn$ is an integer; it is easy to adapt the proof otherwise. Since $h_n \to h_\infty$, for any $\varepsilon>0$ there exists $n_\varepsilon$ such that if $n-j > n_\varepsilon$,
\be
1-\varepsilon < \frac{h_{n-j}}{h_n} < 1 + \varepsilon.
\ee
The last term of \eqref{a useful expanded relation indeed} is then less than $\varepsilon$ if $s<1$ and $n$ large enough (the case $s=1$ is trivial). For any $\delta>0$, there exists $n_\delta$ such that $\sum_{j>n_\delta} \frac1j |1-\e{-\alpha_j}| < \delta$. The sum over the first $n_\delta$ terms in the middle term of \eqref{a useful expanded relation indeed} is then less than $(1+\varepsilon) \frac{n_\delta}n \sum \frac1j |1-\e{-\alpha_j}|$, and it vanishes in the limit $n\to\infty$; the sum over the remaining $n-n_\delta$ terms is less than $(1+\varepsilon) \delta$.
\end{proof}

We conclude this section with a corollary which follows immediately from Proposition \ref{prop estimates h_n} and which will be very useful later.

\begin{corollary}
\label{cor B}
Let
\[
B = \sup_{m,n} \frac{h_m}{h_n}.
\]
Then if $\sum_{\ell\geq1} \frac1\ell |1-\e{-\alpha_\ell}| < \infty$, we have $0<B<\infty$.
\end{corollary}

\section{The spatial model of random permutations}
\label{sec spatial model}

We now introduce the spatial structure and consider the model described in the introduction.
Let $\Lambda$ be an open bounded subset of $\bbR^d$, and let $N$ be the number of ``particles" of the system.
The state space of our model is
\be
\Omega_{\Lambda,N} = \Lambda^N \times \caS_N
\ee
with $\caS_N$ the group of permutations of $N$ elements. $\Omega_{\Lambda,N}$ is equipped with the product of the Borel $\sigma$-algebra on $\Lambda^N$, and the discrete $\sigma$-algebra on $\caS_N$. Let $(\bsx,\pi) \in \Omega_{\Lambda,N}$ with $\bsx = (x_1,\dots,x_N)$; our Hamiltonian is given by \eqref{energy}.

\subsection{The thermodynamic potentials}

We consider the canonical ensemble where the particle density is fixed, and the grand-canonical ensemble where the chemical potential is fixed.
The canonical partition function of this model is
\be
\label{can part fct}
Y(\Lambda,N) = \frac1{N!} \int_{\Lambda^N} \dd\bsx \sum_{\pi \in \caS_N} \e{-H(\bsx,\pi)}.
\ee
The division by $N!$ guarantees that the partition function scales like the exponential of the volume.
The definition makes sense for integer $N$; it is convenient to extend the partition function to noninteger $N$, e.g.\ by linear interpolation.
The grand-canonical partition function is
\be
\label{gd can part fct}
Z(\Lambda,\mu) = \sum_{N\geq0} \e{\mu N} Y(\Lambda,N).
\ee
The parameter $\mu$ is called the chemical potential.
Then we can define the two relevant thermodynamic potentials, the free energy and the pressure:
\be
\label{thermo pots}
\begin{aligned}
&q_\Lambda(\rho) = -\frac1{|\Lambda|} \log Y(\Lambda, |\Lambda|\rho), \\
&p_\Lambda(\mu) = \frac1{|\Lambda|} \log Z(\Lambda,\mu).
\end{aligned}
\ee
The parameter $\rho$ is the density.

We define the ``dispersion relation" $\varepsilon(k)$ by
\be
\e{-\varepsilon(k)} = \int_{\bbR^d} \e{-2\pi\ii kx} \e{-\xi(x)} \dd x.
\ee
For now $\varepsilon(k)$ can be complex; $\varepsilon(0)=0$ and $\Re\varepsilon(k) \geq a|k|^2$ for small $k$. The most relevant case is the Gaussian, $\e{-\xi(x)} = (4\pi\beta)^{-d/2} \e{-|x|^2/4\beta}$. This corresponds to the ideal Bose gas and $\varepsilon(k) = 4\pi^2 \beta |k|^2$. We will suppose in the next subsection that $\varepsilon(k)$ is real.

Note also that the transformation $\alpha_j \mapsto \alpha_j + cj$ translates into
\be
\label{transformation}
\begin{split}
&q_\Lambda(\rho) \mapsto q_\Lambda(\rho) + c\rho; \\
&p_\Lambda(\mu) \mapsto p_\Lambda(\mu-c).
\end{split}
\ee

Next we recall the notion of Fisher convergence \cite{Rue}.
A sequence $(\Lambda_n)$ of domains in $\bbR^d$ converges to $\bbR^d$ in the sense of Fisher if
\begin{itemize}
\item $\lim_{n\to\infty} |\Lambda_n| = \infty$.
\item As $\varepsilon \to 0$,
\[
\sup_n \frac{|\partial_{\varepsilon \diam\Lambda_n} \Lambda_n|}{|\Lambda_n|} \longrightarrow 0,
\]
where $\partial_r \Lambda = \{ x \in \bbR^d : \dist(x,\partial\Lambda) \leq r \}$.
\end{itemize}
This notion is very general.
If $\Lambda$ is bounded with piecewise smooth boundary, then the scaled domains $\Lambda_n = \{ nx : x \in \Lambda \}$ form a Fisher sequence.

For the following three theorems, we always suppose that $\alpha_\ell/\ell$ converges as $\ell \to \infty$. Because of the symmetry \eqref{transformation} we can choose the limit, so we suppose that
\be
\lim_{\ell\to\infty} \frac{\alpha_\ell}\ell = 0.
\ee

The first result is about the infinite volume pressure, that is given by an exact expression.

\begin{theorem}
\label{thm pressure}
For any $\mu \in \bbR\setminus\{0\}$, and any sequence $\Lambda_n$ that converges to $\bbR^d$ in the sense of Fisher, we have
\[
\lim_{n\to\infty} p_{\Lambda_n}(\mu) = p(\mu) \isdefby \sum_{n\geq1} \frac{\e{\mu n - \alpha_n}}n \int_{\bbR^d} \e{-n \varepsilon(k)} \dd k.
\]
Notice that $p(\mu)$ is finite and analytic for $\mu<0$, and that $p(\mu)=\infty$ for $\mu>0$.
\end{theorem}

Next, we have the existence of the thermodynamic limit for the free energy.

\begin{theorem}
\label{thm free energy}
There exists a convex function $q(\rho)$ such that, for any sequence $(\Lambda_n)$ of domains converging to $\bbR^d$ (Fisher), and any sequence $(\rho_n)$ of numbers converging to $\rho\geq0$, we have
\[
\lim_{n\to\infty} q_{\Lambda_n}(\rho_n) = q(\rho).
\]
\end{theorem}

Then $q$ is continuous. It is a standard exercise in analysis to show that the above property is equivalent to uniform convergence of $q_n$ to $q$ on compact intervals. Finally, pressure and free energy are related by Legendre transforms, a property known in statistical mechanics as ``equivalence of ensembles".

\begin{theorem}
\label{thm equivalence ensembles}
The infinite volume pressure and free energy are related as follows:
\[
\begin{split}
&q(\rho) = \sup_\mu \bigl[ \rho\mu - p(\mu) \bigr], \\
&p(\mu) = \sup_\rho \bigl[ \rho\mu - q(\rho) \bigr].
\end{split}
\]
And $q(\rho)$ is analytic except at the critical density
\be \label{critdens}
\rho_{\rm c} = \sum_{n\geq1} \e{-\alpha_n} \int_{\bbR^d} \e{-n\varepsilon(k)} \dd k.
\ee
\end{theorem}

One can check that the critical density is real. It is always finite in dimensions $d \geq 3$; but it may be infinite in $d=1,2$, in which case $q$ is real analytic for all $\rho$ in $[0,\infty)$.

\bfig
\centerline{\includegraphics[width=130mm]{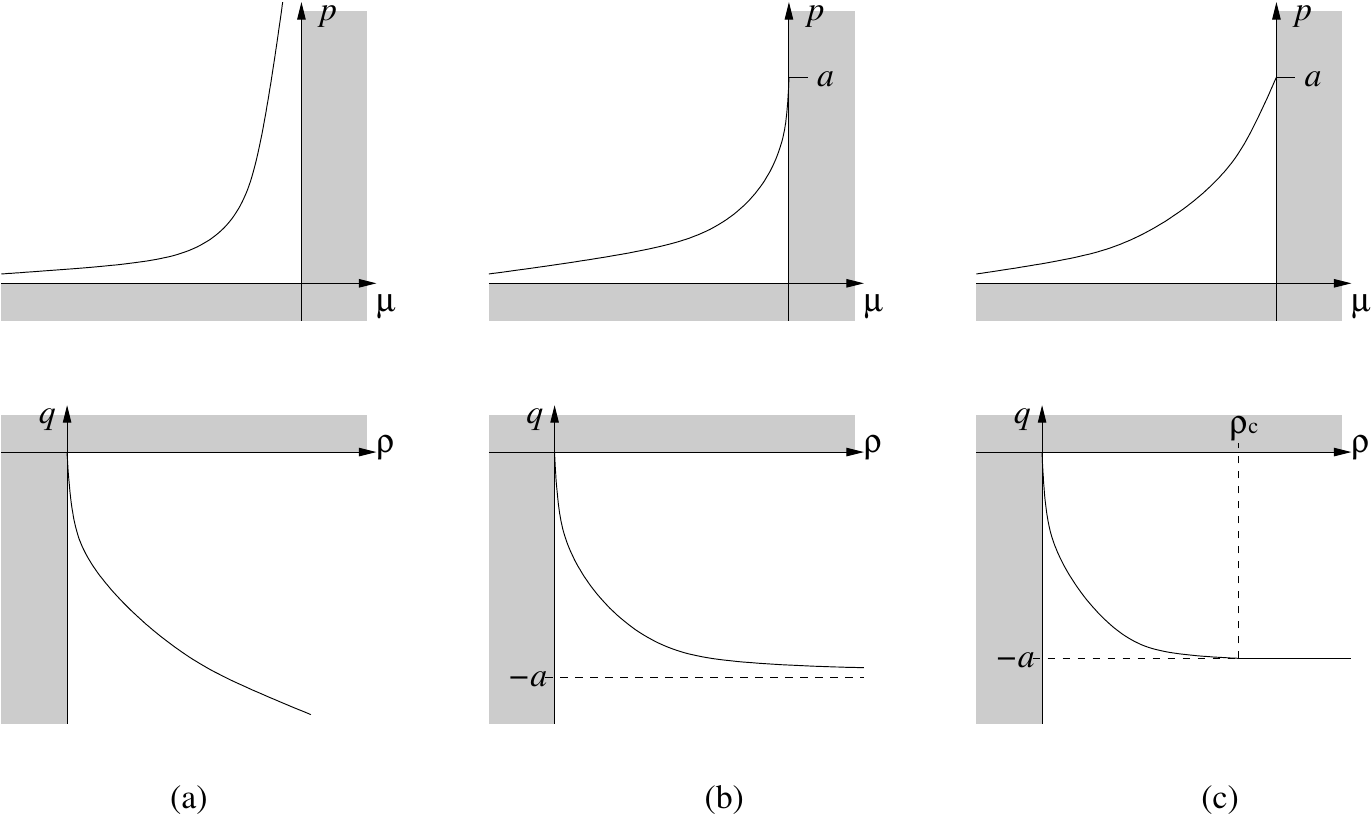}}
\caption{Qualitative graphs of the pressure $p(\mu)$ and its Legendre transform the free energy $q(\rho)$. The critical density $\rho_{\rm c}$ is infinite in (a) and (b), and it is finite in (c).}
\label{fig pq}
\efig

So far we have considered free boundary conditions. For the proofs of our results on the probability model below, we will need  versions of the above theorems with periodized boundary conditions. To be precise, let us consider a sequence $(\Lambda_n)$ of $d$-dimensional boxes with side-length $L_n$, with $L_n \to \infty$. We define $H_\Lambda$ as in (\ref{periodized}) below, and let $q_\Lambda^{\rm per}(\rho)$ and $p_\Lambda^{\rm per}(\mu)$ be the corresponding free energy and pressure. As is usual in statistical mechanics, a change in boundary conditions brings a correction to thermodynamic potentials of the kind $q_\Lambda^{\rm per} - q_\Lambda \approx \frac{|\partial\Lambda|}{|\Lambda|}$, with $|\partial\Lambda|$ a measure of the boundary of $\Lambda$; same for the pressure. The next theorem is less sharp but it is enough for our purpose.

\begin{theorem}
\label{thm periodic bc}
For any $\rho\geq0$ and any $\mu \in \bbR\setminus\{0\}$,
\[
\begin{split}
&\lim_{n\to\infty} q_{\Lambda_n}^{\rm per}(\rho) = q(\rho); \\
&\lim_{n\to\infty} p_{\Lambda_n}^{\rm per}(\mu) = p(\mu).
\end{split}
\]
In both cases, convergence is uniform on compact intervals.
\end{theorem}

While  Theorems \ref{thm pressure}--\ref{thm equivalence ensembles} are proved in Section \ref{sec thermo limit}, Theorem \ref{thm periodic bc} needs the notation of Section \ref{sec Fourier} and it is proved in the appendix.

\subsection{The probability model}
\label{prob model}

We will now study our model from a probabilistic point of view, proving the occurrence of infinite cycles above the critical density. As is often the case, we will need more stringent conditions than we did for studying the thermodynamic potentials.

We also need to slightly modify the Hamiltonian. Let $\Lambda$ be a $d$-dimensional cubic box with side length $L$ and volume $V = L^d$. Define $\xi_{\Lambda}$ through 
\[
\e{-\xi_\Lambda(x)} = \sum_{y \in \bbZ^d} \e{-\xi(x - L y)}. 
\]
The important point is that $\e{-\xi_\Lambda}$ has positive, $\Lambda$-independent Fourier transform $\e{-\eps(k)}$ already in finite volume; this helps us to relate our model with a probability model on Fourier modes. See Proposition \ref{transfer to Fourier} in Section \ref{sec Fourier}. 
In particular, we note that $\int_{\Lambda}{\e{-\xi_\Lambda(x)}} \dd^d x = 1$. If $\e{-\xi}$ has compact support and if $L$ is larger than the diameter of the support, this ``periodized" setting corresponds to usual periodic boundary conditions.

We now define 
\be \label{periodized}
H_\Lambda(\bsx,\pi) = \sum_{i=1}^N \xi_\Lambda(x_i - x_{\pi(i)}) + \sum_{\ell\geq1} \alpha_\ell r_\ell(\pi).
\ee
and introduce a probability measure on $\Omega_{\Lambda,N}$ such that a random variable $\theta: \Omega_{\Lambda,N} \to \bbR$ has expectation
\be
E_{\Lambda,N}(\theta) = \frac1{Y(\Lambda,N) N!} \int_{\Lambda^N} \dd\bsx \sum_{\pi\in\caS_N} \theta(\bsx,\pi) \e{-H_{\Lambda}(\bsx,\pi)}.
\ee
Note that $Y(\Lambda,N)$ is the partition function with periodized boundary conditions, a fact 
that we suppress from the notation. 
 
Next, let $\ell_i(\pi) = 1,2,\dots$ denote the length of the permutation cycle of $\pi$ that contains the index $i$.
It is convenient to consider the density of points that belong to cycles of certain lengths.
Precisely, let
\be
\bsvarrho_{a,b}(\pi) =\frac1{|\Lambda|} \# \{ i = 1,2,\dots : a \leq \ell_i(\pi) \leq b \} = \frac{N_{a,b}(\pi)}{|\Lambda|}.
\ee

\begin{theorem}
\label{thminfinitecycles}
Assume that $\e{-\xi}$ is continuous, that it has positive Fourier transform (i.e.\ $\varepsilon(k)$ is real), and that 
$\sum_{\ell\geq1} \frac{|\alpha_\ell|}{\ell} < \infty$. We also suppose that $\rho_{\rm c} < \infty$. Let $\eta$ be any function such that $\eta(V) \to \infty$ and $\eta(V)/V \to 0$ as $V\to\infty$. Then for all $s\geq0$,
\[
\begin{array}{ll}
\displaystyle \lim_{V\to\infty} E_{\Lambda, \rho V}(\bsvarrho_{1,\eta(V)}) = \begin{cases} \rho & \text{if } \rho \leq \rho_{\rm c}; \\ \rho_{\rm c} & \text{if } \rho \geq \rho_{\rm c}; \end{cases} & \text{(microscopic cycles)} \\[6mm]
\displaystyle \lim_{V\to\infty} E_{\Lambda, \rho V}(\bsvarrho_{\eta(V),V/\eta(V)}) = 0; & \text{(mesoscopic cycles)}\\[4mm]
\displaystyle \lim_{V\to\infty} E_{\Lambda, \rho V}(\bsvarrho_{\eta(V),sV}) = \begin{cases} 0 & \text{if } \rho \leq \rho_{\rm c}; \\ s & \text{if } 0 \leq s \leq \rho-\rho_{\rm c}, \\
\rho - \rho_{\rm c} & \text{if } 0 \leq \rho-\rho_{\rm c} \leq s. \end{cases}
&\text{(macroscopic cycles)} 
\end{array}
\]
\end{theorem}

Recall the expression \eqref{critdens} of the critical density. It is easily seen that
\be
\rho_{\rm c} = \infty \quad \Leftrightarrow \quad \int \frac{\dd k}{\e{\eps(k)}-1} = \infty \quad \Leftrightarrow \quad \int_{|k|<1} \frac{\dd k}{\eps(k)} = \infty.
\ee
The same theorem is stated in \cite{BU} but with $\alpha_\ell \equiv 0$; in this case it holds also when $\rho_{\rm c} = \infty$. The proof of Theorem \ref{thminfinitecycles}, which can be found at the end of Section \ref{sec Fourier}, requires Proposition \ref{prop typical occupation numbers}, which states that certain occupation numbers of Fourier modes are typical. But we can prove one of the claims of Proposition \ref{prop typical occupation numbers} only if $\rho_{\rm c} < \infty$. There is little doubt that all properties also hold true when the critical density is infinite.

\section{Thermodynamic limits and equivalence of ensembles}
\label{sec thermo limit}

In this section we consider the thermodynamic limits of the pressure and of the free energy, and we prove Theorems \ref{thm pressure}--\ref{thm equivalence ensembles}.
We adapt the methods devised in the 1960's by Fisher and Ruelle for classical particle systems.
See \cite{Rue} for references.

\begin{proof}[Proof of Theorem \ref{thm pressure}]
The case $\alpha_2 \neq 0$, $\alpha_n=0$ for any $n\neq2$, was treated in \cite{BU}.
The general case is similar.
The key observation is that everything factorizes according to the permutation cycles.

The grand-canonical partition function reads
\be
Z(\Lambda,\mu) = \sum_{N\geq0} \frac{\e{\mu N}}{N!} \sum_{\pi\in\caS_N} \prod_{n\geq1} \Bigl( \int_{\Lambda^n} \dd x_1 \dots \dd x_n \e{-\sum_{i=1}^n \xi(x_i - x_{i+1})} \e{-\alpha_n} \Bigr)^{r_n(\pi)}
\ee
where $r_n(\pi)$ is the number of cycles of length $n$ in the permutation $\pi$.
We also supposed that $x_{n+1} \equiv x_1$.
The number of permutations of $N$ elements with $r_n$ cycles of length $n$, $n\geq1$, is equal to
\[
N! \Big/ \prod_{n\geq1} n^{r_n} r_n!.
\]
Then
\be
\begin{split}
Z(\Lambda,\mu) &= \prod_{n\geq1} \Bigl\{ \sum_{r\geq0} \frac1{r!} \Bigl( \frac{\e{\mu n - \alpha_n}}n \int_{\Lambda^n} \dd x_1 \dots \dd x_n \prod_{i=1}^n \e{-\xi(x_i - x_{i+1})} \Bigr)^r \Bigr\} \\
&= \exp\Bigl\{ \sum_{n\geq1} \frac{\e{\mu n - \alpha_n}}n \int_{\Lambda^n} \dd x_1 \dots \dd x_n \prod_{i=1}^n \e{-\xi(x_i - x_{i+1})} \Bigr\},
\end{split}
\ee
and thus
\be
\label{expression pressure}
p_{\Lambda}(\mu) = \frac1{|\Lambda|} \sum_{n\geq1} \frac{\e{\mu n - \alpha_n}}n \int_{\Lambda^n} \dd x_1 \dots \dd x_n \prod_{i=1}^n \e{-\xi(x_i - x_{i+1})}.
\ee
Now, with $x_1=0$ in the second term,
\be
\begin{split}
\frac1{|\Lambda|} \int_{\Lambda^n} \dd x_1 \dots \dd x_n \prod_{i=1}^n \e{-\xi(x_i - x_{i+1})} &\leq \int_{\bbR^{d(n-1)}} \dd x_2 \dots \dd x_n \prod_{i=1}^n \e{-\xi(x_i - x_{i+1})} \\
&= \Bigl( \e{-\xi} \Bigr)^{*n}(0) \\
&= \int_{\bbR^d} \e{-n\varepsilon(k)} \dd k.
\end{split}
\ee
The inequality becomes an identity in the limit $\Lambda \nearrow \bbR^d$, and therefore we get Theorem \ref{thm pressure} by dominated convergence if $\mu<0$. Since $\eps(0) = 0$ and $\eps(k) > a |k|^2$, the integrals of $\e{-n\varepsilon(k)}$ around $k=0$ decay as an inverse power of $n$, and the integral for large $k$ decays exponentially. Then if $\mu>0$, we have $p_\Lambda(\mu) = \infty$ for any $\Lambda$ large enough.
\end{proof}

Let $F(\Lambda,N)$ denote the macroscopic free energy,
\be
F(\Lambda,N) = -\log Y(\Lambda,N).
\ee

\begin{lemma}
\label{lem subadditivity}
Let $\Lambda_1, \Lambda_2$ be disjoint open bounded subsets of $\bbR^d$.
Then for any integers $N_1, N_2$,
\[
F(\Lambda_1 \cup \Lambda_2, N_1 + N_2) \leq F(\Lambda_1,N_1) + F(\Lambda_2,N_2).
\]
\end{lemma}

Among the many useful consequences of this subadditive property, we get the upper bound:
\be
\label{upper bound}
F(\Lambda,N) \leq \sum_{i=1}^N F(\Lambda_i,1) = N (\xi(0) + \alpha_1 + \log\rho),
\ee
where the $\Lambda_i$'s partition $\Lambda$ in subdomains, each of volume $|\Lambda|/N = 1/\rho$.

\begin{proof}
We show that $Y(\Lambda_1 \cup \Lambda_2, N_1 + N_2) \geq Y(\Lambda_1,N_1) Y(\Lambda_2,N_2)$.
From the definition,
\be
Y(\Lambda_1 \cup \Lambda_2, N_1 + N_2) = \frac1{(N_1+N_2)!} \int_{(\Lambda_1 \cup \Lambda_2)^{N_1}} \dd\bsx \int_{(\Lambda_1 \cup \Lambda_2)^{N_2}} \dd\bsy \sum_{\pi \in \caS_{N_1+N_2}} \e{-H(\bsx\times\bsy,\pi)}.
\ee
We get an upper bound by restricting the integrals so that exactly $N_1$ points fall in $\Lambda_1$ and $N_2$ points in $\Lambda_2$.
Rearranging the integrals, we obtain
\be
Y(\Lambda_1 \cup \Lambda_2, N_1 + N_2) \geq \frac1{N_1! N_2!} \int_{\Lambda_1^{N_1}} \dd\bsx \int_{\Lambda_2^{N_2}} \dd\bsy \sum_{\pi \in \caS_{N_1+N_2}} \e{-H(\bsx\times\bsy,\pi)}.
\ee
We restrict the sum to permutations of the kind $\pi = \pi_1 \times \pi_2$, with $\pi_1$ a permutation of the first $N_1$ elements and $\pi_2$ a permutation of the last $N_2$ elements.
Then
\be
\label{almost supermultiplicativity}
Y(\Lambda_1 \cup \Lambda_2, N_1 + N_2) \geq \frac1{N_1! N_2!} \int_{\Lambda_1^{N_1}} \dd\bsx \int_{\Lambda_2^{N_2}} \dd\bsy \sumtwo{\pi_1 \in \caS_{N_1}}{\pi_2 \in \caS_{N_2}} \e{-H(\bsx\times\bsy,\pi_1 \times \pi_2)}.
\ee
Finally, we observe that
\be
H(\bsx\times\bsy,\pi_1 \times \pi_2) = H(\bsx,\pi_1) + H(\bsy,\pi_2).
\ee
Integrals in \eqref{almost supermultiplicativity} factorize, yielding the product of partition functions.
\end{proof}

Next we identify the free energy by considering a special sequence of increasing domains.

\begin{lemma}
\label{lem special limit}
Let $C_n$ be the cube of size $2^n$ centered at the origin, and define
\[
q(\rho) = \lim_{n\to\infty} 2^{-dn} F(C_n, \lfloor 2^n \rho \rfloor).
\]
Then
\begin{itemize}
\item[(a)] the above limit exists indeed.
\item[(b)] $q(\rho)$ is convex.
\end{itemize}
\end{lemma}

\begin{proof}
The existence of the limit follows from a standard subadditive argument, where we show that
\be
\label{lower bound for free energy}
\limsup_{n\to\infty} 2^{-dn} F(C_n, \lfloor 2^n \rho \rfloor) \leq \liminf_{n\to\infty} 2^{-dn} F(C_n, \lfloor 2^n \rho \rfloor).
\ee
Fix $k$; we consider $n>2k$.
The cube $C_n$ can be partitioned into $2^{d(n-k)}$ cubes $C_k$.
We put $\lfloor 2^{dk} \rho \rfloor$ points in the first $2^{d(n-k)} - 2^{d(n-2k)}$ cubes, and the remaining $N_{nk}$ points in the remaining $2^{d(n-2k)}$ cubes.
By subadditivity, and the upper bound \eqref{upper bound}, we get
\be
F(C_n, \lfloor 2^{dn} \rho \rfloor) \leq (2^{d(n-k)} - 2^{d(n-2k)}) F(C_k, \lfloor 2^{dk} \rho \rfloor) + N_{nk} \Bigl( \xi(0) + \alpha_1 + \log \frac{N_{nk}}{2^{d(n-k)}} \Bigr).
\ee
The number $N_{nk}$ is not too big:
\be
N_{nk} = \lfloor 2^{dn} \rho \rfloor - (2^{d(n-k)} - 2^{d(n-2k)}) \lfloor 2^{dk} \rho \rfloor \leq 2^{d(n-k)} (1+\rho).
\ee
Then
\be
\limsup_{n\to\infty} 2^{-dn} F(C_n, \lfloor 2^n \rho \rfloor) \leq 2^{-dk} (1 - 2^{-dk}) F(C_k, \lfloor 2^{dk} \rho \rfloor) + 2^{-dk} (1+\rho) (\xi(0) + \alpha_1 + \log(1+\rho)).
\ee
Taking the $\liminf$ $k\to\infty$ in the right side, we obtain \eqref{lower bound for free energy}.

We now prove that $q(\rho)$ satisfies a certain form a continuity, see \eqref{weak continuity} below.
Let $\varepsilon>0$; we study
\be
q(\rho+\varepsilon) = \lim_n 2^{-dn} F(C_n, \lfloor 2^{dn} (\rho+\varepsilon) \rfloor).
\ee
Fix $k$.
There are $2^{d(n-k)}$ cubes $C_k$ in $C_n$.
Put $\lfloor 2^{dk} \rho \rfloor$ points in the first $\lfloor 2^{d(n-k)} (1-\varepsilon) \rfloor$ cubes, and the remaining $N_{nk}$ points in the remaining domain $D_{nk}$.
By subadditivity and the upper bound \eqref{upper bound},
\be
\label{une borne}
F(C_n, \lfloor 2^{dn} (\rho+\varepsilon) \rfloor) \leq \lfloor 2^{d(n-k)} (1-\varepsilon) \rfloor F(C_k, \lfloor 2^{dk} \rho \rfloor) + N_{nk} \bigl( \xi(0) + \alpha_1 + \log \tfrac{N_{nk}}{|D_{nk}|} \bigr).
\ee
We can estimate the last term:
\be
\begin{split}
&N_{nk} = \lfloor 2^{dn} (\rho+\varepsilon) \rfloor - \lfloor 2^{d(n-k)} (1-\varepsilon) \rfloor \lfloor 2^{dk} \rho \rfloor \leq 2^{dn+1} \varepsilon + 2^{d(n-k)} + 2^{dk} \rho, \\
&|D_{nk}| = 2^{dn} - 2^{dk} \lfloor 2^{d(n-k)} (1-\varepsilon) \rfloor \geq 2^{dn} \varepsilon.
\end{split}
\ee
Then
\be
N_{nk} \bigl( \xi(0) + \alpha_1 + \log\tfrac{N_{nk}}{|D_{nk}|} \bigr) \leq 2^{dn} (2\varepsilon + 2^{-k} + 2^{-d(n-k)} \rho) \bigl( \xi(0) + \alpha_1 + \log(2 + \tfrac{2^{-dk}}\varepsilon + \tfrac{2^{-d(n-k)}}\varepsilon ) \bigr).
\ee
We substitute these bounds into \eqref{une borne}, divide by $2^{dn}$, and let $n\to\infty$. We get
\be
q(\rho+\varepsilon) \leq 2^{-dk} (1-\varepsilon) F(C_k, \lfloor 2^{dk} \rho \rfloor) + (2\varepsilon + 2^{-k}) (\xi(0) + \alpha_1 + \log(2 + \tfrac{2^{-dk}}\varepsilon )).
\ee
This is true for any $k$. As $k \to \infty$,
\be
\label{weak continuity}
q(\rho+\varepsilon) \leq (1-\varepsilon) q(\rho) + 2\varepsilon (\xi(0) + \alpha_1 + \log2).
\ee

Finally, we prove a weak version of convexity.
Combined with \eqref{weak continuity}, it implies that $q(\rho)$ is convex indeed.
For any $\rho_1$ and $\rho_2$ such that $2^{dn-1} \rho_i$ is integer for $n$ large enough,
\be
\begin{split}
q(\tfrac12 \rho_1 + \tfrac12 \rho_2) &= \lim_n 2^{-dn} F(C_n, 2^{dn-1} \rho_1 + 2^{dn-1} \rho_2) \\
&\leq \lim_n 2^{-dn} \bigl[ 2^{d-1} F(C_{n-1}, 2^{d(n-1)} \rho_1) + 2^{d-1} F(C_{n-1}, 2^{d(n-1)} \rho_2) \bigr] \\
&= \tfrac12 q(\rho_1) + \tfrac12 q(\rho_2).
\end{split}
\ee
Indeed, we put a density $\rho_1$ of points in half the cubes $C_{n-1}$, and a density $\rho_2$ in the other half, and we used subadditivity.
\end{proof}

We can now prove Theorem \ref{thm free energy} about the convergence of the thermodynamic limit for the free energy.

\begin{proof}[Proof of Theorem \ref{thm free energy}]
First, we show that
\be
\label{thermo limit upper bound}
\limsup_{n\to\infty} \frac1{|\Lambda_n|} F(\Lambda_n,|\Lambda_n|\rho_n) \leq q(\rho).
\ee
Choose a cube $C$, and a number $\rho' > \rho$.
Given $n$, let $M_n$ be the largest integer such that $M_n |C| \rho' \leq |\Lambda_n| \rho_n$.
We pave $\bbR^d$ with translates of the cube $C$.
The volume of cubes inside $\Lambda_n$ is at least $|\Lambda_n| - |\partial_{\diam C} \Lambda_n|$.
Thus for $n$ large enough, the number of cubes inside $\Lambda_n$ is at least 
\be
\frac{|\Lambda_n| - |\partial_{\diam C} \Lambda_n|}{|C|} \geq \frac{|\Lambda_n|}{|C|} \, \frac{\rho_n}{\rho'} \geq M_n.
\ee
By subadditivity of the free energy,
\be
F(\Lambda_n,|\Lambda_n|\rho_n) \leq M_n F(C, \lfloor |C| \rho \rfloor) + F(D,N)
\ee
with
\be
F(D,N) \leq N (\alpha_1 + \log\tfrac N{|D|}).
\ee
Now
\be
\begin{split}
N &= |\Lambda_n|\rho_n - M_n \lfloor |C|\rho \rfloor \leq |\Lambda_n| \rho_n (1 - \tfrac\rho{\rho'} + \tfrac1{|C| \rho'}) + |C|\rho, \\
|D| &= |\Lambda_n| - M_n |C| \geq |\Lambda_n| (1 - \tfrac{\rho_n}{\rho'}).
\end{split}
\ee
Then
\bm
\tfrac1{|\Lambda_n|} F(\Lambda_n,|\Lambda_n|\rho_n) \leq \tfrac{M_n |C|}{|\Lambda_n|} \tfrac1{|C|} F(C, \lfloor |C|\rho \rfloor) \\
+ \Big( \rho_n (1 - \tfrac\rho{\rho'} - \tfrac1{|C|\rho'}) + \tfrac{|C|\rho}{|\Lambda_n|} \Bigr) \Bigl( \alpha_1 + \log \tfrac{\rho_n (1 - \tfrac\rho{\rho'} + \tfrac1{|C| \rho'}) + \tfrac{|C|}{|\Lambda_n|}}{1 - \tfrac\rho{\rho'}} \Bigr).
\end{multline}
Letting $n\to\infty$, the expression simplifies a bit:
\be
\limsup_{n\to\infty} \tfrac1{|\Lambda_n|} F(\Lambda_n,|\Lambda_n|\rho_n) \leq \tfrac\rho{\rho'} \tfrac1{|C|} F(C, \lfloor |C|\rho \rfloor) + \rho (1 - \tfrac\rho{\rho'} - \tfrac1{|C|\rho'}) \Bigl( \alpha_1 + \log \tfrac{\rho (1 - \tfrac\rho{\rho'} + \tfrac1{|C| \rho'})}{1 - \tfrac\rho{\rho'}} \Bigr).
\ee
We can consider the cubes $C_k$ of Lemma \ref{lem special limit} and take the limit $k\to\infty$.
Then we let $\rho' \to \rho$, and we get \eqref{thermo limit upper bound}.

We complete the proof by showing the complementary lower bound, namely
\be
\label{thermo limit lower bound}
\liminf_{n\to\infty} \tfrac1{|\Lambda_n|} F(\Lambda_n,|\Lambda_n|\rho_n) \geq q(\rho).
\ee

A consequence of the limit in the sense of Fisher is that, given $\Lambda_n \nearrow \bbR^d$, there exists $\eta>0$ such that each $\Lambda_n$ contains a translate of a cube of size $\eta \diam\Lambda_n$.
Given $n$, let $k$ be such that $2^k \geq \diam\Lambda_n > 2^{k-1}$.
Then $\Lambda_n$ is contained in a translate of $C_k$, and
\be
\label{volume estimates}
1 \geq \tfrac{|\Lambda_n|}{|C_k|} \geq 2^{-dk} (\eta \, \diam\Lambda_n)^d \geq 2^{-dk} (\eta 2^{k-1})^d = (\tfrac12 \eta)^d.
\ee
By subadditivity, we have
\be
F(C_k, \lfloor 2^{dk} \rho \rfloor) \leq F(\Lambda_n,|\Lambda_n|\rho_n) + F(C_k \setminus \Lambda_n, \lfloor 2^{dk} \rho \rfloor - |\Lambda_n|\rho_n).
\ee
As $n\to\infty$, we also have $k\to\infty$.
Then, taking the $\liminf$ of the expression above, we get
\bm
\label{une borne superieure}
\liminf_{n\to\infty} \tfrac1{|\Lambda_n|} F(\Lambda_n,|\Lambda_n|\rho_n) \geq \liminf_{n\to\infty} \Bigl[ \tfrac{2^{dk}}{|\Lambda_n|} 2^{-dk} F(C_k, \lfloor 2^{dk} \rho \rfloor) \\
- \tfrac{2^{dk} - |\Lambda_n|}{|\Lambda_n|} \tfrac1{2^{dk} - |\Lambda_n|} F(C_k \setminus \Lambda_n, \lfloor 2^{dk} \rho \rfloor - |\Lambda_n|\rho_n) \Bigr].
\end{multline}
Notice that $\lfloor 2^{dk} \rho \rfloor - |\Lambda_n|\rho_n \equiv |C_k \setminus \Lambda_n| \rho_n'$ with
\be
\rho_n' = \tfrac{\lfloor 2^{dk} \rho \rfloor - |\Lambda_n|\rho_n}{|C_k \setminus \Lambda_n|} \longrightarrow \rho
\ee
as $n \to \infty$.
This uses \eqref{volume estimates}.
Using the result for the $\limsup$, Eq.\ \eqref{thermo limit upper bound}, we find that the right side of \eqref{une borne superieure} is larger than $q(\rho)$.
\end{proof}

There remains to prove Theorem \ref{thm equivalence ensembles}.
It is actually enough to show that $p$ is the Legendre transform of $q$ --- since $q$ is convex, it is necessarily the Legendre transform of $p$.

\begin{proof}[Proof of Theorem \ref{thm equivalence ensembles}]
From \eqref{gd can part fct} and \eqref{thermo pots}, we have
\be
\label{finite volume pressure}
p_\Lambda(\mu) = \frac1{|\Lambda|} \log \sum_{\rho \in \bbN/|\Lambda|} \exp\Bigl\{ |\Lambda| \bigl[ \rho\mu - q_\Lambda(\rho) \bigr] \Bigr\}.
\ee
By restricting the sum over a single $\rho$, we get
\be
\label{il suffisait d'y penser}
p_\Lambda(\mu) \geq \rho\mu - q_\Lambda(\rho).
\ee
This holds for any $\Lambda$ and any $\rho$ such that $|\Lambda|\rho$ is integer. It follows that for any $\rho>0$, we can use a suitable Fisher sequence (such as cubic boxes of size $L = \rho^{-1/d} k$, $k \in \bbN$) so as to obtain
\be
p(\mu) \geq \rho\mu - q(\rho).
\ee
This inequality also holds when taking the supremum over $\rho$ in the right side.

The upper bound requires a bit more work. By \eqref{il suffisait d'y penser}, we have
\be
\rho\mu - q_\Lambda(\rho) = \tfrac12 \rho\mu + \tfrac12 \rho\mu - q_\Lambda(\rho) \leq \tfrac12 \rho\mu + p_\Lambda(\tfrac\mu2).
\ee
Let $A$ be a number that is independent of $\Lambda$, to be determined later. Then
\be
p_\Lambda(\mu) \leq \frac1{|\Lambda|} \log \Bigl\{ \sumtwo{\rho \in \bbN/|\Lambda|}{\rho \leq A} \e{|\Lambda| [\rho\mu - q_\Lambda(\rho)]} + \sumtwo{\rho \in \bbN/|\Lambda|}{\rho>A} \e{|\Lambda| [\frac12 \rho \mu + p_\Lambda(\frac\mu2)]} \Bigr\}.
\ee
The latter sum is equal to $\e{|\Lambda| p_\Lambda(\frac\mu2)} \e{\frac12 \mu \lceil A |\Lambda| \rceil} (1 - \e{\frac12 \mu})^{-1}$.

It follows from Theorem \ref{thm free energy} that $q_\Lambda(\rho)$ converges to $q(\rho)$ uniformly on compact intervals. Thus for any $\varepsilon>0$, there exists $\Lambda$ large enough so that
\be
\rho\mu - q_\Lambda(\rho) \leq \rho^* \mu - q(\rho^*) + \varepsilon,
\ee
for any $0 < \rho < A$. Then
\be
p_\Lambda(\mu) \leq \rho^* \mu - q(\rho^*) + \varepsilon + \frac1{|\Lambda|} \log \Bigl\{ A|\Lambda| + (1 - \e{\frac12 \mu})^{-1} \e{|\Lambda| [p_\Lambda(\frac\mu2) - \rho^* \mu + q(\rho^*) - \varepsilon + \frac12 A\mu]}  \Bigr\}.
\ee
We can choose $A$ large enough such that the exponent is negative.
Letting $\Lambda \nearrow \bbR^d$, we see that $p(\mu) \leq \sup_\rho [\rho\mu - q(\rho)]$ indeed.
\end{proof}

\section{The Fourier model of random permutations}
\label{sec Fourier}

We consider now the setting of Section \ref{prob model}. Thus $\Lambda$ is a cube of side length $L$ and volume $V = L^d$. Let $\Lambda^* = \frac1L \bbZ^d$ denote the dual space. The goal of this section is to describe a probability model on $\Omega_{\Lambda,N}^* = (\Lambda^*)^N \times \caS_N$, whose marginal distribution on permutations coincide with the model of spatial permutations. The method is inspired by \cite{Suto2}, who studied the nature of cycles in the ideal Bose gas. The beginning of this section is almost identical to \cite{BU}, the difference being that $\alpha_\ell \neq0$ here.

For a bounded domain $\Lambda$ the new probability space is discrete. We give ourselves a positive continuous function $\varepsilon(k)$, $k \in \bbR^d$. We suppose that $\varepsilon(0)=0$, and that $\varepsilon(k) \geq a|k|^2$ for small $k$. The probability of $(\bsk,\pi) \in \Omega_{\Lambda,N}^*$ is defined by
\be
\label{prob Fourier}
p_{\Lambda,N}(\bsk,\pi) = \begin{cases} \frac1{\widehat Y (\Lambda,N) N!} \e{-\widehat H(\bsk,\pi)} & \text{if } k_i = k_{\pi(i)} \text{ for all } i, \\ 0 & \text{otherwise.} \end{cases}
\ee
Here, the Hamiltonian $\widehat H$ is given by
\be
\widehat H(\bsk,\pi) = \sum_{i=1}^N \varepsilon(k_i) + \sum_{\ell\geq1} \alpha_\ell r_\ell(\pi).
\ee

This model offers an alternate representation to the model of spatial permutations, as far as the permutations are concerned, thanks to the following relation. Let $\delta_{i,j}$ denote Kronecker's delta symbol.

\begin{proposition} \label{transfer to Fourier}
Suppose that $\widehat{\e{-\xi}} = \e{-\varepsilon}$, and define $H_{\Lambda}$ through (\ref{periodized}). 
Then, for any permutation $\pi$,
\[
\int_{\Lambda^N} \e{-H_\Lambda(\bsx,\pi)} \dd\bsx = \sum_{\bsk \in (\Lambda^*)^N} \e{-\widehat H(\bsk,\pi)} \prod_{i=1}^N \delta_{k_i,k_{\pi(i)}}.
\]
\end{proposition}

In particular, it follows that $\widehat Y(\Lambda,N) = Y(\Lambda,N)$, and 
$E_{\Lambda,N}(\Lambda^N \times \{\pi\}) = p_{\Lambda,N}(\Omega_{\Lambda,N}^{\ast},\pi)$ for all permutations $\pi$.

\begin{proof}
We use the following identity, which follows almost directly from well-known relations between Fourier transform and convolution:
\be
\int_{\Lambda^N} \dd\bsx \prod_{i=1}^N \e{-\xi_\Lambda(x_i - x_{\pi(i)})} = \sum_{\bsk \in (\Lambda^*)^N} \prod_{i=1}^N \Bigl[ \delta_{k_i,k_{\pi(i)}} \e{-\varepsilon(k_i)} \Bigr].
\ee
See \cite{BU}, Corollary 5.3, for more details. Notice that the case $N=1$ reduces to Poisson summation formula, and it holds true because $\e{-\xi}$ is continuous. Multiplying both sides by $\e{-\sum \alpha_\ell r_\ell(\pi)}$, one gets the result.
\end{proof}

Next, we introduce occupation numbers. Let $\caN_\Lambda$ be the set of sequences $\bsn = (n_k)$ of integers indexed by $k \in \Lambda^*$, and let $\caN_{\Lambda,N}$ the set of occupation numbers with total number $N$:
\be
\caN_{\Lambda,N} = \Bigl\{ \bsn \in \caN_\Lambda : \sum_{k \in \Lambda^*} n_k = N \Bigr\}.
\ee
To each $\bsk \in (\Lambda^*)^N$ corresponds an element $\bsn \in \caN_{\Lambda,N}$, with $n_k$ counting the number of indices $i$ such that $k_i = k$. Thus we can view $\bsn$ as a subset of $(\Lambda^*)^N$. The probability \eqref{prob Fourier} yields a probability on occupation numbers. Indeed, summing over permutations and over compatible vectors $\bsk$, we have
\be \label{p bsn}
p_{\Lambda,N}(\bsn) = \frac1{Y(\Lambda,N)} \prod_{k \in \Lambda^*} \e{-n_k \varepsilon(k)} h_{n_k},
\ee
with $h_n$ defined in \eqref{def h_n}.

We will obtain some properties on the probability $p_{\Lambda,N}(\bsn)$ below. But we first relate this probability with the lengths of permutation cycles.

\begin{proposition} \label{de-spatialisation}
\[
E_{\Lambda,N}(\bsvarrho_{a,b}) = \frac1V \sum_{\bsn \in \caN_{\Lambda,N}} p_{\Lambda,N}(\bsn) \sum_{k\in\Lambda^*} E_{n_k}(N_{ab}).
\]
\end{proposition}

\begin{proof}
Let $p_{\Lambda,N}(\bsk) = \sum_\pi p_{\Lambda,N}(\bsk,\pi)$. Since $\bsvarrho_{ab}$ depends only on permutations, we have
\be
E_{\Lambda,N}(\bsvarrho_{ab}) = \sum_{\bsk \in (\Lambda^*)^N} p_{\Lambda,N}(\bsk) \sum_{\pi \in \caS_N} \bsvarrho_{ab}(\pi) p_{\Lambda,N}(\pi|\bsk).
\ee
Actually, $\bsvarrho_{ab}(\pi)$ depends only on the conjugacy class of $\pi$. In other words, we have, for any $\pi,\sigma \in \caS_N$,
\be
\bsvarrho_{ab}(\sigma^{-1} \pi \sigma) = \bsvarrho_{ab}(\pi).
\ee
It follows that
\be
\sum_{\pi \in \caS_N} \bsvarrho_{ab}(\pi) p_{\Lambda,N}(\pi|\bsk) = \sum_{\pi \in \caS_N} \bsvarrho_{ab}(\pi) p_{\Lambda,N}(\pi | \sigma(\bsk)).
\ee
Summing first over occupation numbers and then over compatible $\bsk$'s, we get
\be
E_{\Lambda,N}(\bsvarrho_{ab}) = \sum_{\bsn \in \caN_{\Lambda,N}} p_{\Lambda,N}(\bsn) \sum_{\pi \in \caS_N} \bsvarrho_{ab}(\pi) p_{\Lambda,N}(\pi|\bsk).
\ee
Here, $\bsk$ is any vector that is compatible with $\bsn$.

A permutation $\pi$ such that $\pi(\bsk)=\bsk$ (that is, $k_{\pi(i)} = k_i$ for all $i$) can be decomposed into permutations $(\pi_k)_{k \in \Lambda^*}$, where $\pi_k$ is a permutation of the $n_k$ indices $i$ such that $k_i = k$. Notice that
\be
\label{breaking into modes}
N_{ab}(\pi) = \sum_{k \in \Lambda^*} N_{ab}(\pi_k), \qquad r_\ell(\pi) = \sum_{k \in \Lambda^*} r_\ell(\pi_k).
\ee
Then
\be
p_{\Lambda,N}(\bsk,\pi) = \frac1{Y(\Lambda,N) N!} \prod_{k \in \Lambda^*} \e{-n_k \varepsilon(k)} \e{-\sum_\ell \alpha_\ell r_\ell(\pi)}
\ee
if $\pi(\bsk)=\bsk$; it is 0 otherwise. Also,
\be
p_{\Lambda,N}(\bsk) = \frac1{Y(\Lambda,N) N!} \prod_{k \in \Lambda^*} \e{-n_k \varepsilon(k)} h_{n_k} n_k!.
\ee
Then
\be
p_{\Lambda,N}(\pi|\bsk) = \frac{p_{\Lambda,N}(\bsk,\pi) }{p_{\Lambda,N}(\bsk)} = \prod_{k \in \Lambda^*} \frac{\e{-\sum_\ell \alpha_\ell r_\ell(\pi)}}{h_{n_k} n_k!}
\ee
if $\pi(\bsk)=\bsk$, and 0 otherwise. Using \eqref{breaking into modes},
\be
\begin{split}
E_{\Lambda,N}(\bsvarrho_{ab}) &= \frac1V \sum_{\bsn \in \caN_{\Lambda,N}} p_{\Lambda,N}(\bsn) \sum_{(\pi_k \in \caS_{n_k})} \sum_{k \in \Lambda^*} \bsN_{ab}(\pi_k) \prod_{k' \in \Lambda^*} \frac{\e{-\sum_\ell \alpha_\ell r_\ell(\pi_{k'})}}{h_{n_{k'}} n_{k'}!} \\
&= \frac1V \sum_{\bsn \in \caN_{\Lambda,N}} p_{\Lambda,N}(\bsn) \sum_{k \in \Lambda^*} E_{n_k}(N_{ab}).
\end{split}
\ee
\end{proof}

In the light of Proposition \ref{de-spatialisation}, we can now focus on the quantity $p_{\Lambda,N}(\bsn)$. Namely, 
our results from Section \ref{sec cycle weight} imply that macroscopic cycles appear if and 
only if at least one mode is macroscopically occupied, i.e.\ iff 
$p_{\Lambda,N}(n_k \geq s N) > 0$ uniformly in $N \in \N$ and $\Lambda$ such that $N = \rho \Lambda$. 

We prove now that macroscopic occupation can occur only for $k=0$ and that it occurs if and only if $\rho$ is above the critical density $\rho_{\rm c}$ defined in \eqref{critdens}. The first step is a result that gives detailed information about the limiting distribution of the random variable $n_0/V$.

\begin{theorem} \label{exponential RV}
Let $\rho_0 = \max(0,\rho-\rho_{\rm c})$, with $\rho_{\rm c}$ the critical density defined in \eqref{critdens}. Then for all $\lambda \geq 0$, we have
\[
\lim_{V \to \infty} E_{\Lambda, \rho V} ( \e{\lambda n_0 / V} ) = \e{\lambda \rho_0}.
\]
\end{theorem}

\begin{proof}
Our proof is based on the work of Buffet and Pul\'e \cite{BP} for the ideal Bose gas, see also \cite{BU}. We need  to modify it due to the presence of cycle weights. We define 
\be
Y(\Lambda,N,j) = \sum_{\bsn \in \caN_{\Lambda,N}} \biggl( \prod_{k \in \Lambda^\ast} \e{-\eps(k) n_k} h_{n_k} \biggr) \times \begin{cases} h_{n_0 + j}/h_{n_0} & \text{if } j \geq 0, \\ h_{\infty}/h_{n_0} & \text{if } j<0. \end{cases}
\ee
Then $Y(\Lambda,N,0) = Y(\Lambda,N)$.
Recall that $h_\infty = \lim_{n \to \infty} h_n$ is given in Proposition \ref{prop h_n formula} (d) and that $0 < h_\infty < \infty$. Then by \eqref{p bsn}, 
\be \label{n_0 probability}
\begin{split}
p_{\Lambda,N}(n_0 \geq j) &= \frac{1}{Y(\Lambda,N)} \sumtwo{\bsn \in \caN_{\Lambda,N}}{n_0 \geq j} \,\,
\prod_{k \in \Lambda^\ast} \e{- n_k \eps(k)} h_{n_k}\\
&=  \frac{1}{Y(\Lambda,N)} \sum_{\bsn \in \caN_{\Lambda,N-j}} 
\frac{h_{n_0+j}}{h_{n_0}}\prod_{k \in \Lambda^\ast} \e{- n_k \eps(k)} h_{n_k} \\
&= \frac{Y(\Lambda,N-j,j)}{Y(\Lambda,N)}
\end{split}
\ee
for all $j \geq 0$. Using
\be
p_{\Lambda,N}(n_0 = j) = p_{\Lambda,N}(n_0 \geq j) - p_{\Lambda,N}(n_0 \geq j+1),
\ee
the change of summation index $j \mapsto N-j$ gives 
\be \label{fin vol exponential} 
E_{\Lambda,N} (\e{\nu n_0}) = \frac{\e{\nu N}}{Y(\Lambda,N)} \sum_{j=0}^N \e{- \nu j} 
(Y(\Lambda,j,N-j) - Y(\Lambda,j-1,N-j+1)).
\ee
Here we used the convention $Y(\Lambda,-1,N+1) = 0$ and the fact that $p_{\Lambda,N} (n_k \geq N+1) = 0$. We now fix 
$\rho \geq 0$ and put $\nu = \lambda/V$. As a first step we show that 
\be \label{conjecture 1}
\lim_{V \to \infty} E_{\Lambda, \rho V}(\e{\lambda n_0/V}) = 1 \qquad \text{for all } \rho \leq \rho_{\rm c}.
\ee
Above, we wrote $\rho V$ instead of $\floor{\rho V}$, and we will continue to do so in order to simplify the notation. 
To prove \eqref{conjecture 1}, we get from \eqref{fin vol exponential} that 
\be
\label{Gruessgott, Volker!}
\begin{split}
& E_{\Lambda,\rho V}(\e{\frac{\lambda n_0}{V}}) - \e{-\frac{\lambda}{V}} =  \\
& = \frac{\e{\lambda \rho}}{Y(\Lambda,\rho V)} \sum_{j=0}^{\rho V} \e{- \frac{\lambda j}{V}} 
\Big(Y(\Lambda,j,\rho V-j) - Y(\Lambda,j-1,\rho V-j+1)\Big) - \e{-\frac{\lambda }{V}} \\
& = \frac{\e{\lambda \rho}}{Y(\Lambda,\rho V)} \bigg( \sum_{j=0}^{\rho V} \e{- \frac{\lambda j}{V}} 
Y(\Lambda,j,\rho V-j) - \sum_{j=0}^{\rho V} \e{- \frac{\lambda (j+1)}{V}} Y(\Lambda,j,\rho V-j) \biggr)\\
& = \frac{\e{\lambda \rho}}{Y(\Lambda,\rho V)} (1 - \e{-\frac{\lambda}{V}})  \sum_{j=0}^{\rho V} 
\e{- \frac{\lambda j}{V}} 
Y(\Lambda,j,{\rho V}-j)
\end{split}
\ee
We need to show that \eqref{Gruessgott, Volker!} converges to zero as $V \to \infty$. Recall the constant $B$ defined in Corollary \ref{cor B}. Since $Y(\Lambda,N,j) \leq B Y(\Lambda,N)$, we have
\bm
E_{\Lambda,\rho V}(\e{\frac{\lambda n_0}{V}}) - \e{-\frac{\lambda}{V}} \leq B \frac{\e{\lambda \rho}}{Y(\Lambda,\rho V)} \frac{\lambda}{V} \sum_{j=0}^{\rho V} \e{- \frac{\lambda j}{V}} 
Y(\Lambda,j) \\
= B \e{\lambda \rho} \frac{\lambda}{V} \biggl( \sum_{j=0}^{(1-\eps)\rho V} \e{- \lambda j/V} 
 \e{-V (q_{\Lambda}(j/V) - q_{\Lambda}(\rho))} + \sum_{j=0}^{\eps \rho V} \e{- \lambda (\rho V-j)/V} 
 \frac{Y(\Lambda,\rho V-j)}{Y(\Lambda,\rho V)} \biggr)
\end{multline}
for any $\eps > 0$. Above, recall that $q_{\Lambda}(j/V)$ is the finite volume free energy given by \eqref{thermo pots}. It follows from Theorem \ref{thm periodic bc} that $q_{\Lambda}$ converges uniformly on compact intervals to the convex function $q$. By Theorem \ref{thm equivalence ensembles}, $\rho \mapsto q(\rho)$ is strictly decreasing for $\rho < \rho_{\rm c}$. For each $\eps > 0$ there is $b_{\eps} > 0$ such that 
$q_\Lambda(j/V) - q_\Lambda(\rho) > b_\eps$ for all $V$ large enough, and all $j \leq (1-\eps) \rho V$.
So the first term in the bracket above is bounded by $(1 - \eps) \rho V  \e{- b_\eps V}$ and thus converges to zero as $V \to \infty$. For the second term, we claim that $Y(\Lambda, \rho V-j) / Y(\Lambda,\rho V) \leq B$ for all $j$. This is proved by putting the extra $j$ particles into the zero mode $k=0$, or more formally through 
\be \label{Z_N relative bound} 
Y(\Lambda,N) \geq \sumtwo{\bsn \in \caN_{\Lambda,N}}{n_0 \geq j} \prod_{k \in \Lambda^\ast} \e{- \eps(k) n_k} h_{n_k} 
\geq \frac{1}{B} Y(\Lambda,N-j).
\ee
Thus the second term, along with the prefactor $\e{\lambda \rho}\lambda/V$, is bounded by 
$B^2 \lambda \e{\lambda \rho} \rho \eps$. As $\eps$ is arbitrarily small, we obtain \eqref{conjecture 1}.

We now turn to the case $\rho > \rho_{\rm c}$. We define the atomic measure 
\be \label{mu def}
\mu_{\Lambda,\rho} = C_{\Lambda} \sum_{j=0}^{\infty} \Bigl( Y(\Lambda,j,\rho V-j) - Y(\Lambda,j-1,\rho V-j+1) \Bigr) \delta_{j/V}
\ee
on $\bbR^+$, with
\be \label{C_Lambda}
C_\Lambda = \biggl( \sum_{\bsn \in \check\caN_\Lambda} \prod_{k \neq 0} \e{-\eps(k) n_k} h_{n_k} \biggr)^{-1};
\ee
here, we set
\be \label{N check} 
\check\caN_{\Lambda} = \{\bsn \in \caN_{\Lambda}: n_0 = 0\}.
\ee
Later, we will also use the notation $\check\caN_{\Lambda,N} = \{\bsn \in \caN_{\Lambda,N}: n_0 = 0\}$.
As we will see below, $C_\Lambda$ is a correct normalisation so that $\mu_{\Lambda,\rho}$ converges in the limit
$V \to \infty$; it does not depend on $\rho$.

We rewrite \eqref{fin vol exponential} using $\mu_{\Lambda,\rho}$, which gives  
\be \label{integral fraction}
E_{\Lambda,\rho V}(\e{\lambda n_0/V}) = \e{\lambda \rho} \frac{\int 1_{[0,\rho]}(x) \e{-\lambda x} \mu_{\Lambda,\rho}(\dd x)}
{\int 1_{[0,\rho]}(x)  \mu_{\Lambda,\rho}(\dd x)}.
\ee
The strategy is to study the Laplace transform of $\mu_{\Lambda,\rho}$. It will be possible to take the limit $V\to\infty$. Putting the limiting measure in the right side of \eqref{integral fraction}, we will get the infinite volume limit of the left side. An advantage of this strategy is that several convergence issues are handled using standard theorems of analysis. The Laplace transform of $\mu_{\Lambda,\rho}$ is given by 
\be \label{laplace trafo}
\begin{split}
\int_0^\infty \e{-\lambda x} \, \mu_{\Lambda,\rho} (\dd x) &= C_\Lambda (1 - \e{-\lambda/V}) 
\sum_{j=0}^\infty \e{-\lambda j/V} Y(\Lambda,j,\rho V - j) \\
& = C_\Lambda (1 - \e{-\lambda/V}) \sum_{\bsn \in \caN_\Lambda} \biggl( \prod_{k \in \Lambda^\ast} 
\e{- ( \eps(k) + \lambda/V ) n_k} h_{n_k} \biggr) \frac{\tilde h(\bsn,\rho V)}{h_{n_0}},
\end{split}
\ee
with 
\be
\tilde h(\bsn,N) =  
\begin{cases}
h_{(N - \sum_{k \neq 0} n_k)}  & \text{if } \sum_{k \in \Lambda^\ast} n_k \leq N, \\
h_\infty & \text{if } \sum_{k \in \Lambda^\ast} n_k >  N.
\end{cases}
\ee
At this point, the idea in \cite{BP} and \cite{BU} was to factor out the contribution of the zero Fourier mode.
This is not possible here because the factor $\tilde h(\bsn,\rho V)$ couples the modes. This difficulty can be circumvented by introducing
\be \label{nu def}
\nu_{\Lambda} = C_{\Lambda} \sum_{j=0}^{\infty} \Bigl( Y(\Lambda,j,-1) - Y(\Lambda,j-1,-1) \Bigr) \delta_{j/N}
\ee
and by writing
\be \label{measure decomposition}
\int_0^\infty \e{-\lambda x} \mu_{\Lambda,\rho}(\dd x) =  
\int_0^\infty \e{-\lambda x} \nu_{\Lambda}(\dd x) + \int_0^\infty \e{-\lambda x} 
\bigl( \mu_{\Lambda,\rho}(\dd x) - \nu_{\Lambda}(\dd x) \bigr).
\ee
We now prove that the second term vanishes as $V \to \infty$. We have 
\be
\label{some formula}
\begin{split}
& \int_0^\infty \e{-\lambda x} \bigl( \mu_{\Lambda,\rho}(\dd x) - \nu_{\Lambda}(\dd x) \bigr) = C_\Lambda (1 - \e{-\frac{\lambda}{V}}) 
\sum_{j=0}^\infty \e{-\frac{\lambda j}{V}} \Big(Y(\Lambda,j,\rho V - j) - Y(\Lambda,j,-1)\Big) \\
& =  C_\Lambda (1 - \e{-\frac{\lambda}{V}})  \sumtwo{\bsn \in \caN_{\Lambda}}{\sum_k n_k \leq \rho V} 
\biggl( \prod_{k \in \Lambda^\ast} \e{-(\eps(k) + \frac{\lambda}{V}) n_k} h_{n_k} \biggr) 
\frac{\tilde h(\bsn,\rho V) - h_\infty}{h_{n_0}}\\
& =  C_\Lambda (1 - \e{-\frac{\lambda}{V}}) \sum_{n_0 \geq 0} \e{-\frac{\lambda n_0}{V}}   
\sumtwo{\bsn \in \check\caN_{\Lambda}}{\sum_{k \neq 0} n_k \leq \rho V - n_0}
\biggl( \prod_{k \neq 0} \e{-(\eps(k) + \frac{\lambda}{V}) n_k} h_{n_k} \biggr) 
\bigl( \tilde h(\bsn,\rho V) - h_\infty \bigr).
\end{split}
\ee
Now we maximize the second sum in the last line above over $n_0$, which obviously means putting $n_0 = 0$. 
The first sum is a geometric series and cancels the prefactor  $(1 - \e{-\lambda/V})$. As a result, \eqref{some formula} is less than
\[ 
 C_\Lambda    
\sumtwo{\bsn \in \check\caN_{\Lambda}}{\sum_{k\neq 0} n_k \leq \rho V} 
\biggl( \prod_{k \neq 0} \e{-(\eps(k) + \lambda/V) n_k} h_{n_k} \biggr) 
|\tilde h(\bsn, \rho V) - h_\infty|.
\]
The key observation now is that, by Proposition \ref{prop estimates h_n} (a), for given $\eps>0$ there exists $m>0$ such that $|h_{m_0} - h_\infty| < \eps$ whenever $m_0 > m$. Then 
\be
\sumtwo{\bsn \in \check\caN_\Lambda}{\sum_{k \neq 0} n_k \leq \rho V - m} 
\biggl( \prod_{k \neq 0} \e{-(\eps(k) + \lambda/V) n_k} h_{n_k} \biggr) 
|\tilde h(\bsn,\rho V) - h_\infty| \leq 
\eps
\sum_{\bsn \in \check\caN_\Lambda} 
\biggl( \prod_{k \neq 0} \e{-\eps(k) n_k} h_{n_k} \biggr) 
.
\ee
Using the definition of $C_\Lambda$, we find that
\be \label{small piece}
\begin{split}
\int_0^\infty \e{-\lambda x} \bigl| \mu_{\Lambda,\rho}(\dd x) - \nu_{\Lambda}(\dd x) \bigr|  &\leq \eps + C_\Lambda \sum_{N = \rho V - m}^{\rho V}  (h_{\rho V - N} - h_\infty) \check Y_\lambda(\Lambda,N) \\
&\leq 
\eps + B C_\Lambda \sum_{N = \rho V - m}^{\rho V} \check Y_0(\Lambda,N), 
\end{split}
\ee
with 
\be
\check Y_\lambda(\Lambda,N) = \sum_{\bsn \in \check\caN_{\Lambda,N}} \biggl( \prod_{k \neq 0} \e{-(\eps(k) + \lambda/V) n_k} h_{n_k} \biggr).
\ee
Note that 
$
C_\Lambda^{-1} = \sum_{N \geq 1} \check Y_0(\Lambda,N),
$
which suggests that the last term in (\ref{small piece}) is small as the summation from $\rho V - m$ to $\rho V$ 
contains less terms than that giving $C_\Lambda^{-1}$. To prove this, let $\tilde k$ be one of the 
elements of $\Lambda^\ast$ closest to $0$. 
By putting $j$ particles into the mode $\tilde k$, we find that for any $j,N$
\be
\check Y_0(\Lambda,N+j) \geq 
\sum_{\bsn \in \check\caN_{\Lambda,N}} \e{- j \eps(\tilde k) } \frac{h_{n_{\tilde k}+j}}{h_{n_{\tilde k}}} \prod_{k \neq 0} \e{-\eps(k) n_k} h_{n_k} 
\geq B^{-1} \e{- j \eps(\tilde k)} \check Y_0(\Lambda,N).
\ee
Thus for any $N$, 
\be
C_\Lambda^{-1} \geq B^{-1} \check Y_0(\Lambda, N) \sum_{j=0}^\infty \e{- j \eps(\tilde k)} =  B^{-1} \check Y_0(\Lambda, N) 
\frac{1}{1 - \e{-\eps(\tilde k)}} \geq  \frac{\check Y_0(\Lambda, N)}{ B \eps(\tilde k)}.
\ee
Inserting into \eqref{small piece}, we find that
\be
\int_0^\infty \e{-\lambda x} \bigl| \mu_{\Lambda,\rho}(\dd x) - \nu_{\Lambda}(\dd x) \bigr|  \leq \eps + m B^2 \eps(\tilde k).
\ee
As $V\to\infty$ the second term vanishes since $\varepsilon(k)$ is continuous at $0$. Since $\eps$ is arbitrarily small, the left side vanishes indeed in the limit.

Back to \eqref{measure decomposition}. For the first term, 
we can now follow the proof of Theorem A.1 in \cite{BU}. As above, we isolate the contribution of the zero mode and cancel it with the factor $1 - \e{-\lambda/V}$. Thus 
\be
\int_0^\infty \e{- \lambda x} \nu_\Lambda( \dd x) = C_\Lambda h_{\infty} 
\sum_{\bsn \in \check \caN_\Lambda} \biggl( \prod_{k \neq 0} 
\e{- ( \eps(k) + \lambda/V ) n_k} h_{n_k} \biggr).
\ee
Now
\be
\sum_{\bsn \in \check \caN_\Lambda} \biggl( \prod_{k \neq 0} 
\e{- ( \eps(k) + \lambda/V ) n_k} h_{n_k} \biggr) = \exp \biggl( \sum_{k \neq 0} \log \sum_{n \geq 0} 
\e{-(\eps(k) + \frac{\lambda}{V}) n} h_n \biggr).
\ee
By Proposition \ref{prop h_n formula} (e), the logarithm in the exponential is equal to
\be
\sum_{j \geq 1} \e{- ( \eps(k) + \frac{\lambda}{V}) j} 
\frac{\e{-\alpha_j}}{j} = \sum_{j \geq 1} \e{-\eps(k) j} \frac{\e{-\alpha_j}}{j} - 
\frac{1}{V} \int_0^\lambda \sum_{j \geq 1} \e{- ( \eps(k) + \frac{s}{V}) j} \e{-\alpha_j} \, \dd s.
\ee
Using Proposition \ref{prop h_n formula} (e) again, the first term in the inner bracket above is equal to 
$\sum_{k \neq 0} \log \sum_{n=0}^\infty \e{-\eps(k) n} h_n$, and thus it cancels $C_\Lambda$, while the second 
converges to $- \rho_{\rm c}$ as a Riemann sum for every $s \in [0,\lambda]$, when $V \to \infty$. Thus by dominated convergence in $s$, we obtain
\be
\lim_{V \to \infty } \int_0^\infty \e{- \lambda x} \mu_{\Lambda,\rho}( \dd x) = 
\lim_{V \to \infty } \int_0^\infty \e{- \lambda x} \nu_\Lambda( \dd x) = h_{\infty} \e{- \rho_{\rm c}}.
\ee
Thus by the general theory of Laplace transformations, $\mu_{\Lambda,\rho}$ converges to a delta peak of strength 
$h_{\infty}$ at $\rho_{\rm c}$. The claim of the theorem then follows from 
\eqref{integral fraction} for $\rho > \rho_{\rm c}$, and this completes the proof.
\end{proof}

We prove now that the distribution of the random variable $\bsn$ shows typical behaviour. To that end we introduce the three sets
\be
\begin{split}
&A_\epsilon = \bigl\{ \bsn \in \caN_{\Lambda,N} : \bigl| \tfrac{n_0}V - \rho_0 \bigr| < \epsilon \bigr\} \\
&B_{\epsilon,\delta} = \bigl\{ \bsn \in \caN_{\Lambda,N} : \sum_{0<|k|<\delta} n_k < \epsilon V \bigr\} \\
&C_{\epsilon,\delta,M} = \bigl\{ \bsn \in \caN_{\Lambda,N} : \sumtwo{k \in \Lambda^*, |k|\geq\delta}{n_k > M} n_k < \epsilon V \bigr\}.
\end{split}
\ee

\begin{proposition}
\label{prop typical occupation numbers}
Under the assumptions of Theorem \ref{thminfinitecycles}, for any density $\rho$ we have the following.
\begin{itemize}
\item[(a)] For any $\epsilon>0$, $\displaystyle \lim_{V\to\infty} p_{\Lambda, \rho V}(A_\epsilon) = 1$.
\item[(b)] Suppose that $\rho_{\rm c} < \infty$. For any $\epsilon>0$, there exists $\delta_\epsilon$ such that $\displaystyle p_{\Lambda, \rho V}(B_{\epsilon,\delta_\epsilon}) > 1-\epsilon$ for $V$ large enough.
\item[(c)] For any $\epsilon,\delta>0$, there exists $M_{\epsilon,\delta}$ such that $\displaystyle \lim_{V\to\infty} p_{\Lambda, \rho V}(C_{\epsilon,\delta,M_{\epsilon,\delta}}) = 1$.
\end{itemize}
\end{proposition}

The restriction for finite $\rho_{\rm c}$ in item (b) should not be there --- but we cannot prove the claim without it. This is the only reason why Theorem \ref{thminfinitecycles} does not hold when the critical density is infinite.

\begin{proof}
The claim (a) follows from Theorem \ref{exponential RV}:
\be
\lim_{V\to\infty} p_{\Lambda,\rho V}(A_\epsilon) = \lim_{V\to\infty} E_{\Lambda,\rho V} \bigl( 1_{[\rho_0-\epsilon,\rho_0+\epsilon]}(\tfrac{n_0}V) \bigr) = \int 1_{[\rho_0-\epsilon,\rho_0+\epsilon]}(s) \delta_{\rho_0}(s) \dd s = 1.
\ee

We now get a bound on the probability of a given occupation number. Recall the constant $B$ of Corollary \ref{cor B}.
\be
\begin{split}
p_{\Lambda,N}(n_k \geq i) &= \frac1{Y(\Lambda,N)} \sum_{\bsn \in \caN_{\Lambda,N-i}} \e{-\eps(k)i} \biggl( \prod_{k' \in \Lambda^*} \e{-\eps(k') n_{k'}} h_{n_{k'}} \biggr) \frac{h_{n_k+i}}{h_{n_k}} \\
&\leq B \e{-\eps(k)i} \frac{Y(\Lambda,N-i)}{Y(\Lambda,N)} 
\end{split}
\ee
We also have that $\frac{Y(\Lambda,N-i)}{Y(\Lambda,N)} \leq B$, see \eqref{Z_N relative bound}. Then
\be
E_{\Lambda,N}(n_k) = \sum_{i\geq1} p_{\Lambda,N}(n_k \geq i) \leq \frac{B^2}{\e{\eps(k)}-1}.
\ee
By Markov inequality,
\be
p_{\Lambda,N}(B_{\epsilon,\delta_\epsilon}^{\rm c}) \leq \frac{B^2}{\epsilon V} \sum_{0<|k|<\delta_\epsilon} \frac1{\e{\eps(k)}-1} \; \stackrel{V\to\infty}{\longrightarrow} \; \frac{B^2}\epsilon \int_{0<|k|<\delta_\epsilon} \frac{\dd k}{\e{\eps(k)}-1}.
\ee
The integral converges because the critical density \eqref{critdens} is finite. It is possible to choose $V$ large enough and $\delta_\epsilon$ small enough so that $p_{\Lambda,N}(B_{\epsilon,\delta_\epsilon}^{\rm c}) < \epsilon$.

For the claim (c), we use
\be
p_{\Lambda,N}(C_{\epsilon,\delta,M}^{\rm c}) \leq \sum_{m\geq1} \frac1{m!} \sumtwo{k_1,\dots,k_m \in \Lambda^*}{|k_i|>\delta} \frac1{Y(\Lambda,N)} \sumthree{\bsn \in \caN_{\Lambda,N}}{n_{k_i}>M}{\sum_i n_{k_i} \geq \epsilon N} \prod_{k \in \Lambda^*} \e{-\eps(k) n_k} h_{n_k}.
\ee
For given $k_1,\dots,k_m$, we have
\bm
\sumthree{\bsn \in \caN_{\Lambda,N}}{n_{k_i} > M}{\sum n_{k_i} > \epsilon N} \prod_{k\in\Lambda^*} \e{-\eps(k) n_k} h_{n_k} = \\
= \sumtwo{n_1,\dots,n_m > M}{\epsilon N \leq \sum n_i \leq N} \sum_{\bsn' \in \caN_{\Lambda,N-\sum n_i}} \biggl( \prod_{k \in \Lambda^*} \e{-\eps(k) n_k'} h_{n_k'} \biggr) \biggl( \prod_{i=1}^m \e{-\eps(k_i) n_i} \frac{h_{n_{k_i}' + n_i}}{h_{n_{k_i}'}} \biggr).
\end{multline}
We can bound the last ratio by $B$. Then
\be
p_{\Lambda,N}(C_{\epsilon,\delta,M}^{\rm c}) \leq \sum_{m\geq1} \frac{B^m}{m!} \sumtwo{k_1,\dots,k_m \in \Lambda^*}{|k_i|>\delta} \sumtwo{n_1,\dots,n_m > M}{\epsilon N \leq \sum n_i \leq N} \frac{Y(\Lambda, N - \sum n_i)}{Y(\Lambda,N)} \prod_{i=1}^m \e{-\eps(k_i) n_i}.
\ee
We bound the ratio of partition functions by $B$, see \eqref{Z_N relative bound}, and we bound one half of $\eps(k_i)$ by one half of
\be
\eps_0 = \min_{|k|>\delta} \eps(k) > 0.
\ee
Then, since $\sum n_i \geq \epsilon N$, we have
\be
\begin{split}
p_{\Lambda,N}(C_{\epsilon,\delta,M}^{\rm c}) &\leq \e{-\frac12 \eps_0 \epsilon N} \sum_{m\geq1} \frac{B^{m+1}}{m!} \biggl( \sum_{|k|>\delta} \sum_{n>M} \e{-\tfrac12 \eps(k) n} \biggr)^m \\
&\leq B \exp\biggl\{ -V \biggl[ \tfrac12 \eps_0 \epsilon \rho - B \e{-\frac12 \eps_0 M} \frac1V \sum_{|k|>\delta} \frac1{\e{\frac12 \eps(k)}-1} \biggr] \biggr\}.
\end{split}
\ee
We recognise a Riemann sum which is bounded uniformly in $V$. If $M$ is large enough (depending on $\epsilon$ and $\eps_0$, hence on $\delta$), the term in the bracket is positive and everything vanishes in the limit $V\to\infty$.
\end{proof}

We are now equipped for the proof of Theorem \ref{thminfinitecycles}. We use Propositions \ref{de-spatialisation} and \ref{prop typical occupation numbers}, and also Theorem \ref{thm cycle weight}.

\begin{proof}[Proof of Theorem \ref{thminfinitecycles}]
From Proposition \ref{de-spatialisation} we can split
\be
\label{split}
E_{\Lambda,N}(\bsvarrho_{a,b}) = \frac1V \sum_{\bsn \in \caN_{\Lambda,N}} p_{\Lambda,N}(\bsn) \biggl[ E_{n_0}(N_{a,b}) + \sum_{0<|k|<\delta} E_{n_k}(N_{a,b}) + \sum_{|k|\geq\delta} E_{n_k}(N_{a,b}) \biggr].
\ee
We treat the cases separately. By Proposition \ref{prop typical occupation numbers} (a), we can restrict the sum to $\bsn \in A_\epsilon$ with arbitrarily small $\epsilon$. Then for any $\eta(V)$ such that $\eta(V)/V \to 0$,
\be
\lim_{V\to\infty} \sum_{\bsn \in \caN_{\Lambda,\rho V}} p_{\Lambda,\rho V}(\bsn) \tfrac1V E_{n_0}(N_{1,\eta(V)}) = \lim_{V\to\infty} \tfrac1V E_{\rho_0 V}(N_{1,\eta(V)}),
\ee
which is zero by Theorem \ref{thm cycle weight}. On the other hand, for the same reasons we have
\be
\begin{split}
\lim_{V\to\infty} \sum_{\bsn \in \caN_{\Lambda,\rho V}} p_{\Lambda,\rho V}(\bsn) \tfrac1V E_{n_0}(N_{\eta(V),sV}) &= \lim_{V\to\infty} \tfrac1V E_{\rho_0 V}(N_{\eta(V),sV}) \\
&= \begin{cases} s & \text{if } 0 \leq s \leq \rho_0, \\ \rho_0 & \text{if } s \geq \rho_0. \end{cases}
\end{split}
\ee

Next we use Proposition \ref{prop typical occupation numbers} (b) to show that the modes $0<|k|<\delta$ contribute a vanishing amount. Indeed, we can find arbitrarily small $\epsilon$ and $\delta = \delta_\epsilon$ such that
\be
\frac1V \sum_{\bsn \in \caN_{\Lambda,N}} p_{\Lambda,N}(\bsn) \sum_{0<|k|<\delta} E_{n_k}(N_{1,\rho V}) \leq \epsilon + \frac1V \sum_{\bsn \in B_{\epsilon,\delta}} p_{\Lambda,N}(\bsn) \sum_{0<|k|<\delta} E_{n_k}(N_{1,\rho V}) \leq 2\epsilon.
\ee
Therefore we can neglect those modes in \eqref{split} without changing the result in the limit $V\to\infty$.

There remain the modes $|k|>\delta$. By Proposition \ref{prop typical occupation numbers} (c) we can restrict the sum over $\bsn \in C_{\epsilon,\delta,M_{\epsilon,\delta}}$. And because of the definition of $C_{\epsilon,\delta,M_{\epsilon,\delta}}$, we get
\be
\lim_{V\to\infty} \frac1V \sum_{\bsn \in C_{\epsilon,\delta,M_{\epsilon,\delta}}} p_{\Lambda,\rho V}(\bsn) \sum_{|k|>\delta} E_{n_k}(N_{\eta(V),\rho V}) \leq \epsilon.
\ee

The estimates obtained above prove the second and the third claim of Theorem \ref{thminfinitecycles} --- and therefore also the first claim, since the fraction of points in microscopic cycles is obviously equal to the total density, minus the fraction of points in mesoscopic and macroscopic cycles.
\end{proof}

\appendix

\section{Thermodynamic potentials with periodic boundary conditions}
\label{app periodic bc}

We clearly have $Y^{\rm per}(\Lambda,N) \geq Y(\Lambda,N)$, so that
\be
\label{ineq periodic}
q_\Lambda^{\rm per}(\rho) \leq q_\Lambda(\rho).
\ee
It is thus enough to show that $\liminf q_{\Lambda_n}^{\rm per}$ converges to $q$. First we establish some continuity property for $q_\Lambda^{\rm per}$.

\begin{lemma}
\label{lem continuity fen}
For any $\Lambda$ and any $\eta>0$, we have
\[
q_\Lambda^{\rm per}(\rho+\eta) \leq q_\Lambda^{\rm per}(\rho) + \frac{\log B}{|\Lambda|} + \eta \varepsilon_\Lambda(0),
\]
with $B$ the constant of Corollary \ref{cor B}.
\end{lemma}

\begin{proof}
Let $N = \rho |\Lambda|$ and $M = \eta |\Lambda|$. The partition function with periodic boundary conditions is
\be
\begin{split}
Y^{\rm per}(\Lambda,N+M) &= \sum_{\bsn \in \caN_{\Lambda,N+M}} \prod_{k \in \Lambda^*} \e{-n_k \varepsilon_\Lambda(k)} h_{n_k} \\
&\geq \sum_{\bsn \in \caN_{\Lambda,N}} \e{-(n_0+M) \varepsilon_\Lambda(0)} h_{n_0+M} \prod_{k \in \Lambda^* \setminus \{0\}} \e{-n_k \varepsilon_\Lambda(k)} h_{n_k}.
\end{split}
\ee
To get the second line we restricted the sum over occupation numbers to those with $n_0 \geq M$. We know from Corollary \ref{cor B} that $h_{n_0+M} \geq h_{n_0}/B$ and we get the lemma.
\end{proof}

\begin{proof}[Proof of Theorem \ref{thm periodic bc}]
We can derive an expression for the pressure $p^{\rm per}_\Lambda$ like in the proof of Theorem \ref{thm pressure}, namely
\be
p_{\Lambda}^{\rm per}(\mu) = \frac1{|\Lambda|} \sum_{n\geq1} \frac{\e{\mu n - \alpha_n}}n \sum_{k \in \Lambda^*} \e{-n \varepsilon_\Lambda(k)}.
\ee
The last sum is less than $\sum_k \e{-\varepsilon_\Lambda(k)} = \e{-\xi(0)}$. Then $p^{\rm per}_\Lambda$ converges to the expression for $p$ in Theorem \ref{thm pressure} by dominated convergence.

We turn to the free energy with periodic boundary conditions. Let us define
\be
q^{\rm per}(\rho) = \liminf_{n\to\infty} q_{\Lambda_n}^{\rm per}(\rho).
\ee
Suppose that there exists $\rho^*$ such that $q^{\rm per}(\rho^*) < q(\rho^*)$. It follows from Lemma \ref{lem continuity fen} that there exists an interval $I$ close to $\rho^*$ and an $\eta>0$ such that
\be
q^{\rm per}_\Lambda(\rho) < q_\Lambda(\rho) - \eta
\ee
for all $\rho \in I$. Using \eqref{ineq periodic}, we find that the pressure satisfies
\be
p^{\rm per}_\Lambda(\mu) \geq \sum_{N\geq0} \e{\mu N} \e{-q_\Lambda(\frac NV) + \eta 1_I(\frac NV)}.
\ee
Then the infinite volume limit of $p_\Lambda^{\rm per}(\mu)$ is larger than the Legendre transform of $q - \eta \upchi_I$, hence larger than $p(\mu)$ for some $\mu<0$. This contradicts the first claim of Theorem \ref{thm periodic bc}. This shows that $q_\Lambda^{\rm per}$ converges pointwise to $q$. The uniform convergence on compact sets follows from \eqref{ineq periodic}, the uniform convergence of $q_\Lambda$, and Lemma \ref{lem continuity fen}.
\end{proof}

\medskip
{\footnotesize
{\bf Acknowledgments:} We are grateful to the referees for useful comments.
D.U. is grateful to the hospitality of the Erwin Schr\"odinger Insitut of Vienna, the University of Geneva, ETH Z\"urich, and the Center of Theoretical Studies of Prague, where parts of this project were carried forward. V.B.\ is supported by the EPSRC fellowship EP/D07181X/1 and D.U. is supported in part by the grant DMS-0601075 of the US National Science Foundation.
}

\end{document}